\documentclass[reqno]{amsart}
\usepackage{epsfig,amsmath,amsfonts,latexsym,appendix}
\usepackage{amsthm}
\usepackage[abs]{overpic}
\usepackage[usenames,dvipsnames]{xcolor}
\usepackage{palatino}
\usepackage[hyperfootnotes=false]{hyperref}
\usepackage{caption}
\usepackage{dsfont}
\usepackage{mathtools}
\usepackage{cite}
\usepackage{tikz-cd}
\usepackage{amssymb}
\usepackage{comment}
\usepackage{stackrel}
\usepackage{tcolorbox}

\newtheorem{fed}{\textbf{Definition}}[section]
\newtheorem{thm}[fed]{\textbf{Theorem}}

\newtheorem{rem}[fed]{\textbf{Remark}}

\usepackage{amssymb,graphicx,epsfig,psfrag,epic,eepic,latexsym,color}
\usepackage{amsmath}
\usepackage{mathrsfs}
\usepackage{dsfont}

\title{Symplectic methods in the numerical search of orbits in real-life planetary systems\\ \tiny (with numerical explorations of the Jupiter-Europa and Saturn-Enceladus systems)}
\author{Urs Frauenfelder, Dayung Koh, Agustin Moreno}

\address[Urs Frauenfelder]{Institut f\"ur Mathematik \\ Augsburg Universität \\ Augsburg  \\ Germany}

\email{urs.frauenfelder@math.uni-augsburg.de}

\address[Dayung Koh]{Pasadena, California, USA}

\email{dayung.koh@gmail.com}

\address[Agustin Moreno]{School of Mathematics \\ Institute for Advanced Study \\  Princeton  \\ USA \\ \newline Mathematisches Institut \\ Universität Heidelberg \\ Heidelberg \\ Germany}

\email{agustin.moreno2191@gmail.com}

\date{}

\begin{document}

\maketitle

\begin{abstract}
    The intention of this article is to illustrate the use of methods from symplectic geometry for practical purposes. Our intended audience is scientists interested in orbits of Hamiltonian systems (e.g.\ the three-body problem). The main directions pursued in this article are: (1) given two periodic orbits, decide when they can be connected by a regular family of such; (2) use numerical invariants from Floer theory which help predict the existence of orbits in the presence of a bifurcation; (3) attach a sign $\pm$ to each elliptic or hyperbolic Floquet multiplier of a closed \emph{symmetric} orbit, which generalizes the classical Krein--Moser sign to also include the hyperbolic case; and (4) do all of the above in a visual, easily implementable and resource-efficient way. The mathematical framework is provided by the first and third authors in \cite{FM}, where, as it turns out, the ``Broucke stability diagram'' \cite{Broucke} was rediscovered, but further refined with the above signs, and algebraically reformulated in terms of GIT quotients of the symplectic group. If two orbits lie in two different regions of the diagram, they \emph{cannot} be joined by a \emph{regular} family, i.e.\ one should expect a bifurcation. If two orbits have different signs, they also \emph{cannot} be connected by a regular family, even if they lie in the same region of the diagram. The advantage of the above framework is that it applies to the study of closed orbits of an \emph{arbitrary} Hamiltonian system. Moreover, in the case where the system admits symmetries in the form of ``reflections'', i.e.\ anti-symplectic involutions, which is the case for many systems of interest, the information provided for orbits which are symmetric is richer, and one may distinguish more symmetric orbits. This is the case for several well-known families in the space mission design industry, such as the Halo orbits, which are ubiquitous in real-life space missions. We will carry out numerical work based on the \emph{cell-mapping method} as described in \cite{KABM}, for the Jupiter-Europa and the Saturn-Enceladus systems. These are currently systems of interest, falling in the agenda of space agencies like NASA, as these icy moons are considered candidates for harbouring conditions suitable for extraterrestrial life.
\end{abstract}

\tableofcontents

\section{Introduction}

The study of closed orbits of Hamiltonian systems and their bifurcations, is one of the central topics of Floer theory, as introduced by Andreas Floer in a series of papers \cite{F1,F2,F3,F4,F5,F6}; and of Symplectic Field Theory (SFT), as proposed by Eliashberg--Givental--Hofer in \cite{EGH}. Formidable in their depth and scope, both theories underlie many of the powerful methods of modern symplectic geometry. On the other hand, the search of orbits entails significant practical interest. For instance, the restricted three-body problem (3BP), concerning the gravitational motion of a negligible mass around two larger masses, is a prototypical problem in astronomy and is relevant for space mission design. In this context, the influence on a satellite of a planet which comes with an orbiting moon can be approximated by a three-body problem of restricted type. Finding families of orbits for placing such a satellite around the target moon, which minimizes orbit corrections and risk of collisions, is then of central importance for space exploration. The need of organizing all information pertaining to known families of orbits naturally leads to the realm of Big Data, where the numerous modern methods of data analysis (e.g.\ machine learning) apply, and for which computationally cheap methods are highly relevant. Our general direction is then encapsulated in the following guiding questions:

\begin{tcolorbox}
    \textbf{(Classification)} Can we tell when two orbits are \textbf{qualitatively} different\footnote{We say that two orbits are \emph{qualitatively} different if there is no regular family joining them.}? 
    
    \textbf{(Catalogue/Data Science)} Can we resource-efficiently refine data bases of known orbits, and use techniques from data science to study them (e.g.\ machine learning)? 
   
    \textbf{(Practical tests)} Can we use Floer-theoretical invariants to test the accuracy of the algorithms, and to guide/organise the numerical work?

\end{tcolorbox}

The first two questions were addressed in \cite{FM}, where the mathematical framework was set out. In this article, we include numerical work, and we address the third question. We will combine various tools, including:
\begin{enumerate}
    
    \item[(1)] \textbf{Floer numerical invariants:} Euler characteristics of suitable Floer homology groups (one for general closed orbits, and another one which applies for \emph{symmetric} closed orbits); 
    \item[(2)] \textbf{The B-signature} \cite{FM}: a generalization of the classical Moser--Krein signature \cite{kre1,kre2,K3,K4,Moser}, which originally applies only to elliptic Floquet multipliers\footnote{Recall that the Floquet multipliers of a closed orbit are by definition the eigenvalues of the monodromy matrix.}, to also include the case of hyperbolic multipliers, whenever the corresponding orbit is \emph{symmetric}. 
    \item[(3)] \textbf{Global topological methods:} the \emph{GIT-sequence} \cite{FM}, which refines Broucke's stability diagram \cite{Broucke} by adding the $B$-signature.
\end{enumerate}

This paper is the outgrowth of an interdisciplinary dialogue,
whose theme centers on whether methods from modern symplectic geometry (e.g.\ Floer homology) can be of help for engineering problems. While Floer homology was designed to prove statements about the existence of periodic orbits for large classes of Hamiltonian systems, a mere existence statement is of little interest for engineering. However, the situation changes when instead of looking at global Floer homology one looks at \emph{local} Floer homology, as it remains invariant under bifurcations of orbits. If the orbits found so far do not satisfy this invariance requirement then one is sure that there are more families and it makes sense to invest man and computer power to actually detect them. While computing the full local Floer homology might be hard, we instead focus on numerical invariants which are extracted from the full homology and which are easy to implement. This is the motivation to consider the \emph{SFT-Euler characteristic}, as the Euler characteristic of local Floer homology, which can be computed with only knowledge of the spectrum of the monodromy matrix. In the case where the system admits symmetries in the form of involutions, one can consider a further numerical invariant for symmetric orbits, the \emph{real Euler characteristic}, which also stays invariant under a bifurcation. We remark that many orbits which have been found via numerical explorations of classical problems are actually symmetric. As we shall explain, the combination of these \emph{Floer numerical invariants} together with the $B$-signature, provides good tools in order to decide whether to look for periodic orbits, as well as gives useful hints concerning where to actually look for them.

From the side of applications, this dialogue was initiated by the need of engineers of having a firm mathematical groundwork in order to study bifurcations of orbits, for the purpose of space mission design, in connection with real-life planetary systems. Such is the case of the Jupiter-Europa or the Saturn-Enceladus systems, as these icy moons are believed to be candidates for harboring extraterrestrial life, which makes them of tremendous current interest for space agencies such as NASA. We will therefore carry out numerical work for these systems. While we have not attempted to find new orbits in this article, the numerical results in this paper match the predictions of the mathematical framework. 

\smallskip

\textbf{Acknowledgements.} A.\ Moreno is supported by the National Science Foundation under Grant No.\ DMS-1926686, and by the Sonderforschungsbereich TRR 191 Symplectic Structures in Geometry, Algebra and Dynamics, funded by the DFG (Projektnummer 281071066 – TRR 191). The authors are grateful to the referee for very useful comments and corrections on an earlier version of this article.

\section{Preliminaries}

\textbf{Monodromy matrix and symmetric orbits.} Recall that the monodromy matrix of a closed orbit $x$, given by $D_{x(0)}\Psi_T$ where $T$ is the period of $x$ and $\Psi_t$ is the flow, is a symplectic matrix\footnote{A symplectic matrix is a matrix $M$ such that $M^tJM=J$, where $J$ is the standard rotation $J=\left(\begin{array}{cc}
    0 & I \\
   -I  & 0
\end{array}\right)$.}. If the system of  $n+1$ degrees of freedom is described by a \emph{time-independent} Hamiltonian, the eigenvalue $1$ will appear twice in the monodromy matrix (one corresponding to the direction of the flow, another for fixing the energy). We will ignore such trivial eigenvalues, and consider the remaining ones, i.e.\ we consider the \emph{reduced} monodromy matrix, obtained by restricting the dynamics to a level set of the Hamiltonian, and forgetting the direction of the flow\footnote{In practice, it is admittedly simpler to work directly with the unreduced version, and simply work with the eigenvalues. We choose the reduced version for the purpose of exposition. Everything that follows can be easily adapted to the nonreduced case.}. 

Consider a Hamiltonian $H:\mathbb R^{2n+2}\rightarrow \mathbb R$, defined on the phase-space $\mathbb{R}^{2n+2}$ of position-momentum pairs $(q,p)\in \mathbb{R}^{2n+2}$, which comes with the standard symplectic form $\omega=\sum_jdq_j\wedge dp_j$. Consider also an \emph{antisymplectic involution} $\rho$, i.e.\ a map of $\mathbb{R}^{2n+2}$ satisfying $\rho^2=\mathrm{id}, \rho^*\omega=-\omega$. An example of such a map is $$\rho: \mathbb{R}^6\rightarrow \mathbb{R}^6,\;(q_1,q_2,q_3,p_1,p_2,p_3) \mapsto (q_1,-q_2,-q_3,-p_1,p_2,p_3).$$ Assume that $H$ is invariant under $\rho$, i.e.\ $H \circ \rho=H.$ A prototypical example is the Hamiltonian for the restricted three-body problem; see Section \ref{sec:numerics} below. A periodic orbit $x: S^1\rightarrow \mathbb{R}^{2n+2}$ of $H$ is \emph{symmetric} if it satisfies $x(t)=\rho(x(-t)), t \in S^1,$ so that in particular $x\big(0\big),\,\,x\big(\tfrac{1}{2}\big) \in L:=\mathrm{Fix}(\rho)=\{(q,p)\in \mathbb{R}^{2n}:\rho(q,p)=(q,p)\}$ lie in the fixed-point set of $\rho$, and we call them the \emph{symmetric} points. 

\begin{figure}
    \centering
    \includegraphics{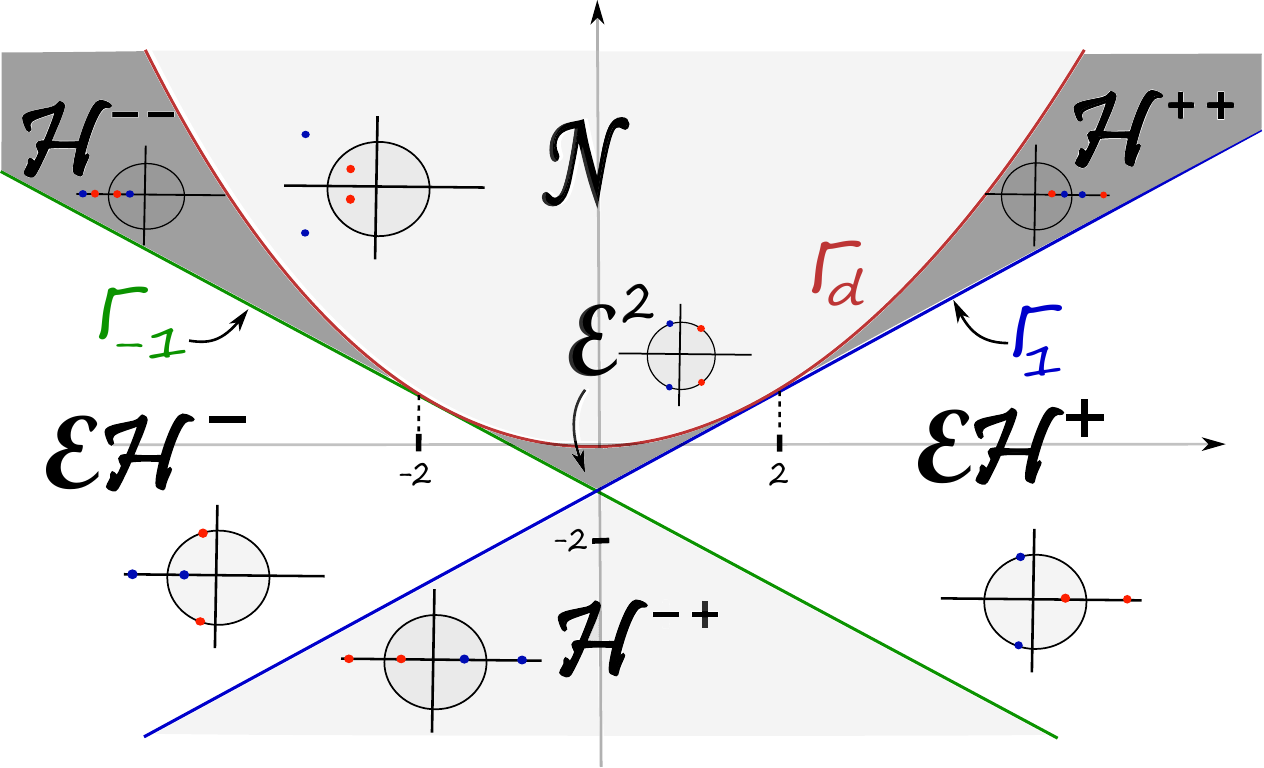}
    \caption{The picture shows $\mathbb{R}^2$, the base of the GIT sequence for $n=2$. This diagram is also referred to as the ``Broucke stability diagram'' \cite{Broucke}. Each region corresponds to a configuration of eigenvalues of the monodromy matrix. We depict a generic configuration for each region; see also Figure \ref{fig:complete_bifurcation}.}
    \label{fig:basen2}
\end{figure}

The reduced monodromy matrix along a symmetric point is a special type of $2n\times 2n$ symplectic matrix $M$. If a basis of $\mathbb{R}^{2n}$ is chosen so that the reflection is the standard one $$\rho=\left(\begin{array}{cc}
    I & 0 \\
    0 & -I
\end{array}\right) \in M_{2n\times 2n}(\mathbb{R}),$$ where $M_{2n\times 2n}(\mathbb{R})$ denotes the space of $2n\times 2n$ real matrices, and $I \in M_{n\times n}(\mathbb{R})$ is the identity matrix, then in this basis $M$ is of the form
\begin{equation}\label{symsymp}
M=M_{A,B,C}=\left(\begin{array}{cc}
A & B\\
C & A^T
\end{array}\right)\in M_{2n\times 2n}(\mathbb{R}).
\end{equation}
Here, $A,B,C$ are $n\times n$-matrices that satisfy the equations
\begin{equation}\label{eq}
B=B^T,\quad C=C^T,\quad AB=BA^T,\quad
A^TC=CA,\quad A^2-BC=I,
\end{equation}
ensuring that $M$ is symplectic. The expression for $M_{A,B,C}$ implies the choice of a basis for the tangent space to the fixed-point locus along the symmetric point. A different choice of basis amounts to acting with an invertible matrix $R \in GL_n(\mathbb{R})$, via
\begin{equation}\label{act}
R_*\big(A,B,C\big)=\Big(RAR^{-1},RBR^T,(R^T)^{-1}CR^{-1}\Big),
\end{equation}
i.e.\ $M_{A,B,C}$ is replaced by $M_{R_*(A,B,C)}$. Moreover, the eigenvalues of $M$ are completely determined by those of the first block $A$, i.e.\ one can reduce the characteristic polynomial of $M$ to that of $A$ (a fact which already appears in \cite{Broucke}).\footnote{Concretely, the characteristic polynomial of $M_{A,B,C}$ is given by
$$p_{A,B,C}(t)=t^n p_{-2A}\big(-t-\tfrac{1}{t}\big),$$
where $p_{-2A}$ is the characteristic polynomial of the matrix $-2A$.} For $n=2$, the distinct cases for eigenvalues of $A$ are summarized in Figure \ref{fig:basen2}, which shows the plane $\mathbb{R}^2$ parametrizing the point $p=(\mathrm{tr}(A),\det(A))\in \mathbb{R}^2$, the coefficients of the characteristic polynomial of $A$. This is explained in more detail in \cite{FM}, where moreover preferred normal forms for each type of matrix is provided (see also \cite{Broucke}). The purpose of the $B$-signature is to refine the data provided by the point $p$, as follows.

\medskip

\textbf{B-signature for symmetric orbits.} We now explain the notion of $B$-signature, introduced in \cite{FM}. For the case $n=2$, the following construction will assign a pair $\epsilon=(\epsilon_1,\epsilon_2)$, the \emph{$B$-signature}, where $\epsilon_i=\pm$ is a ``plus'' or a ``minus'' label for $i=1,2$, to the eigenvalues of the $2\times 2$-block $A$ of $M=M_{A,B,C}$, assuming that they are real and different. We do not assign a sign to the remaining cases, i.e.\ when they coincide or are complex, which correspond respectively to $p$ lying in $\Gamma_d$ or $\mathcal{N}$ in Figure \ref{fig:basen2}. If $\mu_1<\mu_2$ are these eigenvalues, let $v_i$ be an eigenvector of $A^T$ with eigenvalue $\mu_i$, i.e.\ $A^Tv_i=\mu_iv_i$. The \emph{$B$-sign} $\epsilon_i$ of $\mu_i$ is then defined as $$\epsilon_i=\mbox{sign}(v_i^TBv_i)\in \{\pm \},$$ where we use the $B$-block of $M$. It is easily seen that $\epsilon_i$ is independent of the choice of $v_i$.

In the presence of a period doubling bifurcation of a symmetric orbit, the $B$-signature jumps either at $t=0$ or at $t=1/2$ (as in Figure \ref{fig:sym_bif}). The orbit arising from a period doubling bifurcations will be symmetric near the point where the $B$-sign did \emph{not} jump. This will be explicitly checked in  numerical examples below, and serves as a hint as to where to expect bifurcations.


\medskip

\textbf{Floer numerical invariants.} The Floer numerical invariants are, simply put, numbers which stay invariant before and after a bifurcation, when considering families of orbits. They are a count of the number of orbits with suitable signs depending on the Floquet multipliers of each orbit. The invariance of this number implies that they can be used as a practical test for the algorithm used: if this number agrees before and after a bifurcation (as it should), one can rest assured. If the numbers do not agree, one knows that there is at least one orbit missing. The fact that these numbers are invariant is non-trivial, and follows from results in Floer theory. However, if one accepts this as a fact, then this can be used in practice. Of course, one needs a simple way of computing such numbers. We will provide a more detailed guideline in Section \ref{sec:SFT-Euler}. The formula for computing these numbers depends on the dimension, as this determines the number of Floquet multipliers. For example, in dimension $6$, the dimension relevant for the spatial three-body problem, the invariant of a periodic orbit $x$ of a \emph{time-independent} Hamiltonian system, which we call the \emph{SFT-Euler characteristic}, is

\begin{eqnarray*}
\chi_{SFT}(x)&=\#\big\{\mathcal H^{--}\,,\mathcal{EH}^-\,,\mathcal E^2\,,\textrm{ good}\,\,\mathcal H^{++}\,,\mathcal N\big\}\\
&-\#\big\{\mathcal H^{-+},\,\textrm{good}\,\,\mathcal{EH}^+\big\}.
\end{eqnarray*}

Here, $\mathcal{H}^{\pm},\mathcal{E},\mathcal{N}$ stand for \emph{positive/negative hyperbolic}, \emph{elliptic}, and \emph{non-real} orbit. We have also ignored the eigenvalue $1$ (which appears twice), and considered the remaining eigenvalues (two pairs of them, which explains the notation). The ``good'' $\mathcal{EH}^+/\mathcal{H}^{++}$ orbits are those which respectively are \emph{not} even covers of $\mathcal{EH}^-/\mathcal{H}^{-+}$ ones. The $4$-dimensional case is simply given by
\begin{eqnarray*}
\chi_{SFT}(x)=\#\big\{\textrm{ good}\,\,\mathcal H^{+}\big\}-\#\big\{\mathcal E,\,\mathcal{H}^-\big\},
\end{eqnarray*}
where the good $\mathcal{H}^+$ orbits are those which are not even covers of $\mathcal{H}^-$ ones. In short, the invariant is easily computable from the knowledge of the multipliers. In the case where the orbit is symmetric, we will also consider the \emph{real Euler characteristic}. This invariant serves as a further test, as it can detect when the orbit bifurcates as a \emph{chord}, even when it does not bifurcate as an orbit. Its definition is left for Section \ref{sec:real_Euler}.


\medskip

\textbf{Global topological methods.} In practice, one would like to compare the orbits at hand with those that have been found before. The natural notion of equivalence of two orbits is to say that they are \emph{qualitatively} the same, provided one can find a path of orbits joining them, which is regular, i.e.\ it does not undergo bifurcation (the parameter for the path is usually the energy or e.g.\ a mass parameter). For this purpose, we will illustrate the use of global topological methods, via the \emph{GIT sequence} introduced in \cite{FM} by the first and third authors. This is a sequence of three spaces (``top'', ``middle'' and ``base'') consisting of equivalence classes of symplectic matrices, and concrete maps between them, which are also explicitly computable. A closed orbit of an \emph{arbitrary} Hamiltonian system induces a point in the ``base'' and ``middle'' spaces of this sequence. The ``base'' (for $n=2$) is a copy of the plane $\mathbb{R}^2$, split into components labelled according to the Floquet multipliers of the orbit; see Figure \ref{fig:basen2}. The resulting diagram, as we learned after rediscovering it in the context of the GIT sequence, was originally introduced by Broucke in \cite{Broucke} (see also Howard--Mackay \cite{HM} for higher-dimensional versions). 

\begin{figure}
    \centering
    \includegraphics[width=0.7\linewidth]{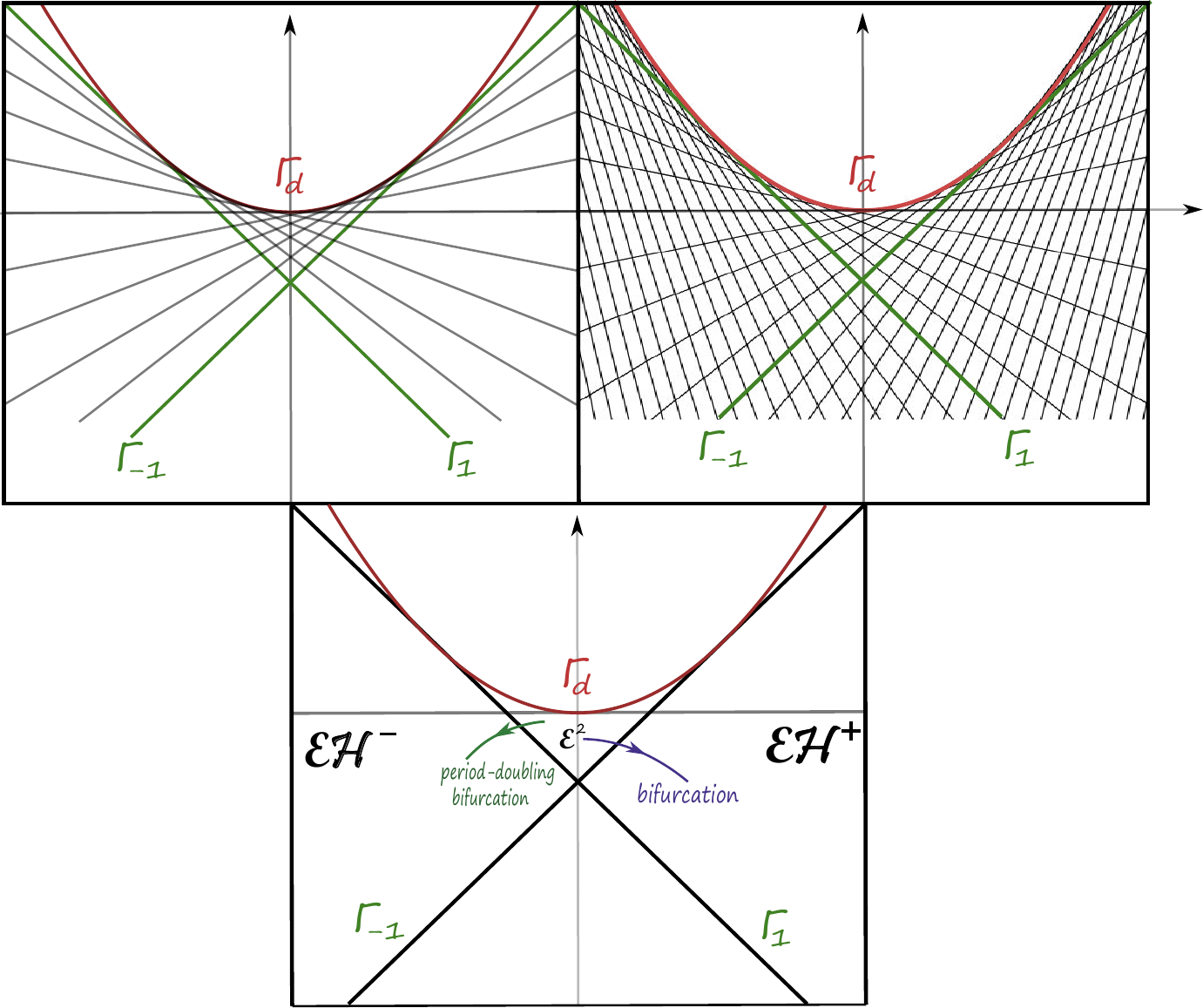}
    \caption{The picture shows $\mathbb{R}^2$, the base of the GIT sequence for $n=2$. The parabola $\Gamma_d=\{y=1/4x^2\}$ corresponds to matrices with double eigenvalues. The locus of matrices with a fixed eigenvalue is a line tangent to $\Gamma_d$. On the left, we have the \emph{elliptic} pencil of lines $\{\Gamma_\theta, \theta \in [0,2\pi)\}$, where $\Gamma_\theta$ has slope $\cos(2\pi \theta)$ and corresponds to the eigenvalue $e^{2\pi i\theta}$. On the right, the complete pencil, also containing the \emph{hyperbolic} pencil $\{\Gamma_\lambda: \lambda\in \mathbb R\backslash[-1,1]\}$, where $\Gamma_\lambda$ has slope $a(\lambda)=\frac{1}{2}(\lambda+\frac{1}{\lambda})$ and corresponds to hyperbolic eigenvalue $\lambda$. A family of orbits bifurcating induces a family of monodromy matrices whose eigenvalues crosses $1$, and hence is seen as a path in $\mathbb{R}^2$ which crosses $\Gamma_1$; we sketch such a path in the bottom panel, where one eigenvalue pair goes from elliptic to positive hyperbolic, while the other stays elliptic. We also sketch an example of a period-doubling bifurcation, where one eigenvalue pair goes from elliptic to negative hyperbolic. Similarly, a $k$-fold bifurcation, where $\lambda=e^{2\pi i \frac{l}{k}}$ (with $\lambda^k=1$) is being crossed, induces a path crossing $\Gamma_{l/k}$, for some $l$. We will plot numerical examples below.}
    \label{fig:complete_bifurcation}
\end{figure}

If the orbit is symmetric, and we choose one of its symmetric points, then there is an associated point in the ``top'' space. A family of orbits induces a path in the corresponding space of the sequence. These spaces also contain subsets corresponding to bifurcations of orbits (which look like a pencil of lines tangent to a parabola; see Figure \ref{fig:complete_bifurcation}). A family of orbits which bifurcates induces a path in the GIT spaces which crosses the component of the bifurcation loci corresponding to the type of bifurcation. Therefore if two orbits correspond to points that lie in different regular components (each of the seven of Figure (\ref{fig:basen2})), there is no regular family that joins them. In other words, the topology of these spaces can be used as an obstruction to the existence of regular families of orbits, as well as the study of bifurcations, in a concrete, visual manner. What is more, one can refine the obstructions provided by the bifurcation loci, by attaching the $B$-signature to the point corresponding to an orbit. Even if two orbits induce points lying in the same regular component of the plane, there is no regular family between them if they have different pair of sign labels. See Appendix \ref{app:GIT} for a mathematical treatment of the GIT sequence.  

\medskip

\textbf{Data base.} The topological approach also serves the purpose of providing a \emph{data base} for orbits, in the form of a cloud of dots in the plane with labels attached. The ``data point'' that is then associated to a symmetric orbit with monodromy matrix $M_{A,B,C}$ is the tuple $(p=(\mathrm{tr}(A),\det(A)),\epsilon=(\epsilon_1,\epsilon_2))$, which is independent of all choices, except perhaps the choice of symmetric point at which we linearize (here, $\epsilon$ is empty for the complex or double-eigenvalue case). For a given orbit, the data stored is relatively cheap, which makes the approach resource-efficient. 

\medskip

\medskip

\textbf{Non-symmetric orbits and Krein theory.} In practice, it may not be apparent whether a given orbit is symmetric. However, in order to refine the data proved by the point $p$, we can still appeal to classical Krein theory, as proposed by Krein \cite{kre1,kre2}\cite{K3,K4} and rediscovered by M\"oser \cite{Moser}, which associates a sign only to the \emph{elliptic} eigenvalues. It turns out that, in the elliptic case, the notion of $B$-signature coincides with the more classical one of Krein signature; see \cite{FM}. 
In the case the orbit is not symmetric, one may still compute the associated point $p=(\mathrm{tr}(A),\det(A))$ via the expressions 
$$
\det (A)=\frac{a}{4}-\frac{1}{2},\;\mathrm{tr}(A)=\frac{b}{2},
$$
where, if $\lambda_1^p,\lambda_2^p,\lambda_1^s,\lambda_2^s$ are the eigenvalues of $M$ (satisfying $\lambda_1^s\lambda_2^s=\lambda_1^p\lambda_2^p=1$), then $a,b$ are given by the symmetric polynomials
$$
a=\lambda^p_1\lambda_1^s+\lambda^p_1\lambda_2^s+\lambda^p_2\lambda_1^s+\lambda^p_2\lambda^s_2+2,
$$
and
$$
b=\lambda^p_1+\lambda^p_2+\lambda^s_1+\lambda^s_2=\mathrm{tr}(M).
$$
Note that these expressions can be computed directly from the eigenvalues of $M$ without need of conjugating $M$ to be of the form $M_{A,B,C}$. If one or both pairs of the eigenvalues are elliptic (and distinct), i.e.\ $p$ lies in a region of Figure \ref{fig:basen2} with at least one $\mathcal E$ label on it, we attach a sign to each of them as follows. We first complexify, i.e.\ we now work over $\mathbb C$. Consider the matrix $$G=-iJ=\left(\begin{array}{cc} 0 & -i Id\\ i Id& 0\end{array}\right)$$ which we think of as acting on $\mathbb C^4$ by left-multiplication (each block is $2\times 2$). Let $\lambda=e^{i\theta},\overline{\lambda}=e^{-i\theta}$ be an elliptic pair of eigenvalues of $M$, with im$(\lambda)>0$. Given $v\in \mathbb{C}^4$ an eigenvector for $\lambda$, the \emph{Krein-sign} of $\lambda$ is 
$$
\kappa(\lambda)=\mbox{sign}(v^tG\overline v).
$$ It is easy to check that $v^tG\overline v$ is a non-zero real number, whose sign is independent of $v$.\footnote{Indeed, since $\lambda$ is simple, any other vector is of the form $w=\mu\cdot v$ with $\mu\neq 0$, and so $w^tG\overline w=\vert \mu \vert^2 v^tG\overline v$.} We have the property $\kappa(\overline{\lambda})=-\kappa(\lambda)$.

\medskip


\section{The SFT-Euler characteristic of a periodic orbit}\label{sec:SFT-Euler}

The SFT-Euler characterstic of a closed orbit can be defined in any dimension, as the Euler characteristic of its so-called \emph{local Floer homology} (which we will not treat here). If the orbit is degenerate and \emph{good} (i.e.\ \emph{not} an even multiple cover of a negative hyperbolic orbit), then this invariant is $\pm 1$, depending on the parity of its so-called \emph{Conley-Zehnder} index; if it is non-degenerate and \emph{bad}, it is zero (we shall not need the definition of the Conley-Zehnder index, but the interested reader can consult e.g.\ \cite{FvK}). The interesting case is when it is degenerate, and hence may undergo bifurcation. If one adds a small perturbation, it might bifurcate into a finite collection of other non-degenerate orbits. One then defines the SFT-Euler characterstic as the number of good orbits with even Conley-Zehnder index, minus the number of good orbits with odd Conley-Zehnder index, where the orbits are taken among the orbits that appear after perturbation. The remarkable fact is that this number is independent on the perturbation. In particular, it remains invariant before and after a bifurcation. In what follows, we shall explain this in more concrete terms, in dimensions $4$ and $6$.

\subsection{The four-dimensional case}

In the four-dimensional case, the reduced monodromy matrix $M$ is a $2\times 2$-matrix. Since it is
symplectic its determinant is $1$ and therefore its spectrum is completely determined by its trace. The trace of the reduced monodromy matrix can as well be obtained from the trace of the nonreduded monodromy matrix $\widetilde{M}$ by the formula
$$\mathrm{tr}(M)=\mathrm{tr}(\widetilde{M})-2.$$
We distinguish the following cases.
\begin{description}
 \item[Negative hyperbolic case] In this case 
 $$\mathrm{tr}(M)\leq -2.$$
 If $\mathrm{tr}(M)<-2$ the spectrum of $M$ contains two negative real eigenvalues which are inverse
 to each other. If $\mathrm{tr}(M)=-2$, then the only eigenvalue of $M$ is $-1$, having algebraic
 multiplicity two.  
 \item[Elliptic case] In this case 
 $$-2<\mathrm{tr}(M)<2.$$
 The spectrum of $M$ contains two nonreal eigenvalues on the unit complex circle which are inverse to each other.
 \item[Degenerate case] In this case 
 $$\mathrm{tr}(M)=2.$$
 The spectrum of $M$ only contains $1$.
 \item[Positive hyperbolic case] In this case 
 $$\mathrm{tr}(M)>2.$$
 The spectrum of $M$ contains two positive real eigenvalues different from one which are inverse to each other. 
\end{description}
We can now explain the parity of the Conley-Zehnder index, which is roughly speaking a rotation number associated to $M$, and its relationship to good/bad orbits. In the negative hyperbolic as well as
the elliptic case the parity of the Conley-Zehnder index is odd, while in the positive hyperbolic case, it
is even. In the degenerate case bifurcation occurs, and we do not define the parity of the Conley-Zehnder index. To define the SFT-Euler characteristic we additionally need the distinction of periodic orbits into good and bad ones. In the four-dimensional case, periodic orbits whose Conley-Zehnder index has
odd parity are always \emph{good}, i.e.\ negative hyperbolic and elliptic orbits are always good. On the other
hand, positive hyperbolic ones can be bad. To explain what this means we need to recall that a periodic orbit gives rise to multiple covers of itself. The monodromy matrix of the $k$-fold cover is then $M^k$ and if $\lambda$ is an eigenvalue of $M$, then $\lambda^k$ is an eigenvalue of $M^k$. In particular, 
if $M^k$ is negative hyperbolic or elliptic, the same is true for $M$. On the other hand, nondegenerate
even covers of negative hyperbolic ones are positive hyperbolic. If a positive hyperbolic orbit is an even cover of a negative hyperbolic one it is called \emph{bad}. Therefore, in dimension four, the
SFT-Euler characteristic of an orbit $x$ is defined by
\begin{eqnarray*}
\chi_{\mathrm{SFT}}(x)&= \#\big\{\textrm{good orbits with even CZ-index}\big\} \\
&-\#\big\{\textrm{orbits with odd CZ-index}\big\}\\
&=\#\big\{\textrm{good positive hyperbolic orbits}\big\}\\
&-\#\big\{\textrm{elliptic and negative hyperbolic orbits}\big\}.
\end{eqnarray*}
That the SFT-Euler characteristic before and after a bifurcation does not change follows from
the invariance of local Floer homology. In Appendix \ref{app:invariance}, we explicitly check this for generic and some non-generic examples. 

\subsection{The six-dimensional case}

In the six-dimensional case the reduced monodromy matrix $M$ is a $4\times 4$-matrix. 
Since it is symplectic its characteristic polynomial 
$$p(x)=x^4+c_3 x^3+c_2 x^2+c_1 x+1$$
is a palindrome, i.e.\ $c_1=c_3$. In particular, there exists a quadratic matrix
$$q(y)=y^2+b_1 y+b_0$$
such that
$$p(x)=x^2 q\big(x+\tfrac{1}{x}\big).$$
Let $\mu_1$ and $\mu_2$ be the roots of the quadratic polynomial $q$. Since the polynomial
$q$ is a real polynomial its roots are both real or complex conjugate to each other. We now distinguish
several cases with the help of the roots of $q$.
\begin{description}
 \item[Nonreal case ($\mathcal N$)] The two roots $\mu_1$ and $\mu_2$ are not real. In this case the eigenvalues of
 the monodromy matrix $M$ are neither real nor lie on the unit circle. They appear as a quadruple
 $\big(\lambda, \tfrac{1}{\lambda}, \overline{\lambda}, \tfrac{1}{\overline{\lambda}}\big)$. 
\end{description}
The real case has to be subdivided into several subcases. If $\mu_1$ and $\mu_2$ are real and distinct
then maybe after a symplectic change of coordinates the reduced monodromy matrix splits as
$$M=\left(\begin{array}{cc}
M_1 & 0\\
0 & M_2
\end{array}\right)$$
where $M_1$ and $M_2$ are symplectic $2\times 2$-matrices satisfying
$$\mathrm{tr}(M_1)=\mu_1, \quad \mathrm{tr}(M_2)=\mu_2.$$
The roots $\mu_1$ and $\mu_2$ hence determine if $M_1$ respectively $M_2$ is elliptic, negative hyperbolic, positive hyperbolic or degenerate. In the real case we order the roots such that 
$$\mu_1 \leq \mu_2.$$
\begin{description}
 \item[Doubly negative hyperbolic case $\mathcal H^{--}$] In this case we have
 $$\mu_1 \leq -2, \quad \mu_2 \leq -2.$$
 \item[Elliptic/ negative hyperbolic case $\mathcal{EH}^-$ ]In this case we have
 $$\mu_1 \leq -2, \quad -2<\mu_2<2.$$
 \item[Negative/ positive hyperbolic case $\mathcal H^{-+}$ ] In this case we have
 $$\mu_1 \leq -2, \quad \mu_2>2.$$
 \item[Doubly elliptic case $\mathcal E^2$ ]In this case we have
 $$-2 <\mu_1<2,  \quad -2<\mu_2<2.$$
 \item[Elliptic/ positive hyperbolic case $\mathcal{EH}^+$ ] In this case we have
 $$-2<\mu_1<2, \quad \mu_2>2.$$
 \item[Doubly positive hyperbolic case $\mathcal H^{++}$ ] In this case we have
 $$\mu_1>2, \quad \mu_2>2.$$
 \item[Degenerate case $\mathcal D$ ] In this case we have $\mu_1=2$ or $\mu_2=2$. 
\end{description}
The parity of the Conley-Zehnder index is additive. For example since the parity of the Conley-Zehnder index in the elliptic case is odd and in the positive hyperbolic case is even it follows that in the
$\mathcal{EH}^+$-case it is odd again. The following table displays the parity in the various cases.
\\ \\
\begin{center}
\begin{tabular}{|c|c|}
\hline
Type & Parity of Conley-Zehnder index\\
\hline\
$\mathcal H^{--}$ & even\\
\hline\
$\mathcal{EH}^-$ & even\\
\hline\
$\mathcal H^{-+}$ & odd\\
\hline\
$\mathcal E^2$ & even\\
\hline \
$\mathcal{EH}^+$ & odd\\
\hline \
$\mathcal H^{++}$ & even\\
\hline \
$\mathcal N$ & even\\
\hline
\end{tabular}
\end{center}
\smallskip
The only periodic orbits which can be bad are the ones of $\mathcal{EH}^+$ and $\mathcal H^{++}$ type. Namely a periodic
orbit of $\mathcal{EH}^+$ type is bad, if it is an even cover of one of $\mathcal{EH}^{-}$ type. Similarly, a periodic orbit
of $\mathcal H^{++}$-type is bad, if it is an even cover of an $\mathcal H^{-+}$-orbit. Otherwise orbits are good. For example if an $\mathcal H^{++}$-orbit is an even cover of a $\mathcal H^{--}$ orbit it is good. The SFT-Euler characteristic of the orbit $x$ is now
\begin{eqnarray*}
\chi_{SFT}(x)&=  \#\big\{\textrm{good orbits with even CZ-index}\big\} \\
&-\#\big\{\textrm{good orbits with odd CZ-index}\big\}\\
&=\#\big\{\mathcal H^{--}\,,\mathcal{EH}^-\,,\mathcal E^2\,,\textrm{ good}\,\,\mathcal H^{++}\,,\mathcal N\big\}\\
&-\#\big\{\mathcal H^{-+},\,\textrm{good}\,\,\mathcal{EH}^+\big\}.
\end{eqnarray*}

\section{The real Euler characteristic of a symmetric orbit}\label{sec:real_Euler}

In this section, we focus on the particular case of symmetric orbits, i.e.\ orbits which are invariant under an antisymplectic involution which preserves the Hamiltonian, as defined in the Introduction. We will follow the exposition of \cite{FvK2}, where the Hörmander index is introduced.

A symmetric periodic $x$ can be seen both as a periodic orbit, as well as a chord between the Lagrangian fixed-point locus of the involution. Therefore it has a Conley-Zehnder index $\mu_{CZ}(x)$, and a \emph{Lagrangian Maslov index} $\mu_L(x)$, which is a \emph{half}-integer, i.e.\ takes values in $\frac{1}{2}\mathbb{Z}$ (again, its definition will not be needed, but can be found e.g.\ in \cite{RS}). The difference of these two, as introduced in \cite{FvK2}, is the \emph{Hörmander index} $$s(x) = \mu_{CZ}(x) - \mu_L(x)\in \frac{1}{2}\mathbb{Z},$$ also a half-integer. One can use this index to detect when $x$ bifurcates as a chord, even when it doesn't bifurcate as an orbit; we will explain this in the next section via concrete examples.

We note that the iterates of a symmetric periodic orbit $x^k$ for $k \in \mathbb N$ are symmetric orbits as well. We say that a symmetric periodic orbit is nondegenerate if for any $k \in \mathbb N$ we have $\det(\Phi^k - I) \neq 0$, where $\Phi$ is the monodromy matrix of $x$, i.e.\ $1$ is not an eigenvalue of any iterate of $\Phi$. Moreover, the Chebyshev polynomials of the first kind are recursively defined by
$$
T_0(x) = 1
$$
$$
T_1(x) = x
$$
$$
T_{k+1}(x) = 2xT_k(x) - T_{k-1}(x).
$$
The Chebyshev polynomials of the second kind are similarly defined by
$$
U_0(x) = 1
$$
$$
U_1(x) = 2x
$$
$$
U_{k+1}(x) = 2xU_k(x) - U_{k-1}(x).
$$
The following gives a formula for computing the Hörmander index of the iterates of a symmetric orbit, in terms of the monodromy matrix, which in particular is easy to implement numerically, and does not make use of the definition of the Conley-Zehnder index nor the Lagrangian Maslov index.
\begin{thm}\cite{FvK2}\label{thm:sindex} Let $x$ be a nondegenerate, symmetric periodic orbit, with monodromy matrix
$$
M=M_{A,B,C}=\left(\begin{array}{cc}
A & B\\
C & A^T
\end{array}\right),
$$
satisfying Equations (\ref{eq}) given in the Introduction. Then the Hörmander indices of its iterates are given by
\begin{equation}\label{eq:sindex}
s(x^k)= \frac{1}{2}\mbox{sign}\left((Id -T_k(A))U_{k-1}(A)^{-1}C^{-1}\right),
\end{equation}
$k \in \mathbb N$. For $k=1$, we have in particular that
\begin{equation}\label{eq:sindex1}
s(x) = \frac{1}{2}\mbox{sign}\left((Id-A)C
^{-1}\right).
\end{equation}
\end{thm}
Here, $sign$ denotes the signature of a matrix (the number of positive eigenvalues, minus the number of negative eigenvalues). We have also used the fact that $C$ is invertible if $x$ is non-degenerate \cite[Lemma 3.2]{FvK2}. 

\smallskip

\textbf{Real Euler characteristic.} Similarly to the Conley-Zehnder index (which induces the SFT-Euler characteristic), one can consider the Euler characteristic of the so-called \emph{local Lagrangian Floer homology} of a symmetric orbit $x$, when viewed as a chord. We call the resulting quantity the \emph{real Euler characteristic} $\chi_L(x)$. More concretely, this works as follows. Before or after a bifurcation, one obtains a collection of non-degenerate symmetric orbits, for which computes the parity of the Maslov index $\mu_L(x)$. By this, if $\mu_L(x)=\frac{1}{2}m_L(x)$, we mean the parity of $m_L(x)\in \mathbb Z$; note that $m_L(x)$ is even if and only if $\mu_L(x)$ is an integer. Note that, in practice, without needing to know the definition of this index, its parity can be determined from the following: 
\begin{itemize}
    \item The monodromy matrix;
    \item the formula $\mu_L(x)=\mu_{CZ}(x)-s(x)$;
    \item Formula (\ref{eq:sindex}) (or (\ref{eq:sindex1})) of Theorem \ref{thm:sindex}, which in particular gives the parity of $s$;
    \item the table above giving the parity of the CZ-index (which is always an integer) in terms of the eigenvalue classification of the monodromy matrix.
\end{itemize}

The real Euler characteristic $\chi_L(x)$ is then defined as
$$
\chi_L(x)=\sum_j(-1)^{\mu_L(x_j)}=\sum_i(-1)^{\mu_{CZ}(x_j)}(-1)^{-s(x_j)}=\sum_ii^{m_L(x_j)} \in \mathbb{C},
$$
where the sum runs over the collection $x_j$ of non-degenerate chords arising after perturbation of $x$. Note that by definition, $\chi_L(x)$ is complex-valued. Its invariance under bifurcation follows from invariance of the local Lagrangian Floer homology of $x$.

\section{symmetric period-doubling}\label{sec:examples} In this section we discuss period-doubling of symmetric orbits, in dimension four, corresponding to the case where two eigenvalues collide at $-1$. This will illustrate the use of the invariants and Krein-type signatures that we have discussed. While the discussion is general, we will illustrate it in concrete numerical examples in Section \ref{sec:numerics}. 

\smallskip

\begin{figure}
    \centering
    \includegraphics[width=1 \linewidth]{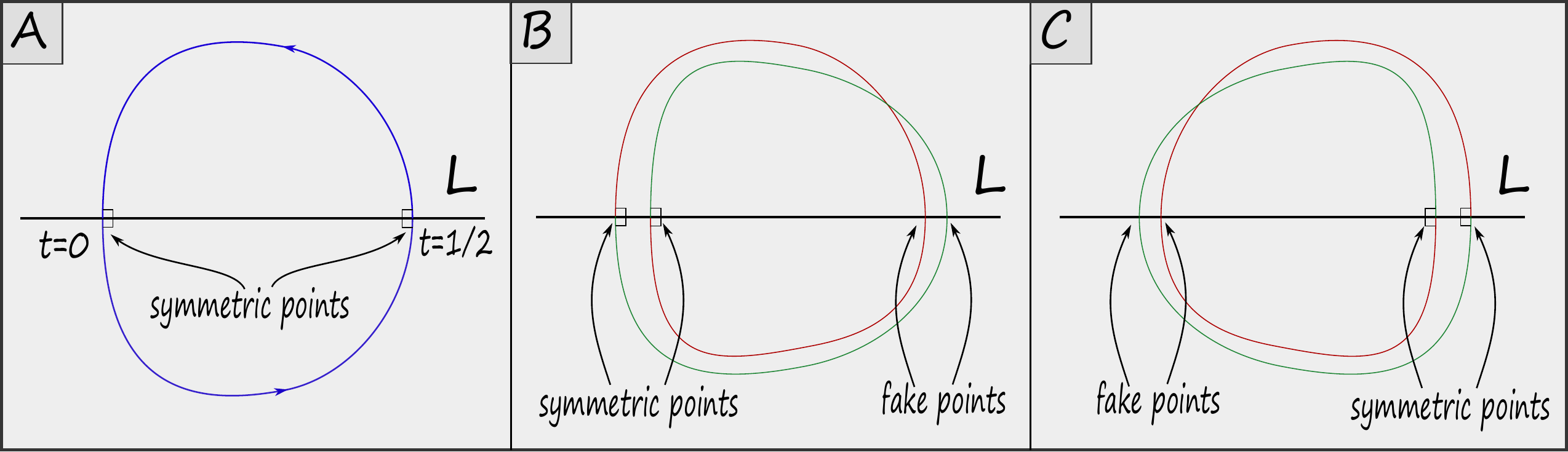}
    \caption{(A) A simple symmetric orbit; we indicate the intersection with $L=$Fix$(\rho)$ (i.e.\ the symmetric points) with right angles. (B) A period doubling bifurcation, only intersecting $L$ near $x(0)$ (we also indicate the ``fake'' intersection points near $x(1/2)$, at which the new orbit \emph{misses} $L$). (C) The alternative version of (B), intersecting $L$ only near $x(1/2)$.}
    \label{fig:sym_bif}
\end{figure}

\textbf{Symmetric period-doubling.} We study the case corresponding to \cite[p.\ 599]{abraham-marsden}, but in the symmetric case. Consider a simple symmetric periodic orbit $x$ which intersects the Lagrangian fixed-point locus at time $t=0$ and $t=1/2$, and which belongs to a family whose reduced monodromy matrix goes from elliptic to negative hyperbolic; see Figure \ref{fig:sym_bif}(A). As a simple orbit there is no bifurcation since there is no eigenvalue $1$ in
the reduced monodromy matrix. However, we can interpret this orbit as a chord from the Lagrangian to itself, where  $x(0)$ happens to agree with $x(1)$; now, as a \emph{chord}, it might bifurcate.\footnote{Such a chord bifurcation takes place if and only if the Lagragian Maslov index jumps, which happens if and only if the Hörmander index jumps. See Section \ref{sec:real_Euler} for definitions.} If this happens, we can apply the symmetry again to the red chord in Figure \ref{fig:sym_bif}(B)/(C) to obtain the green chord in the same figure. So, \emph{two} chords bifurcate. This is compatible with the real Euler characteristic. Indeed, the Lagrangian Maslov index of $x$ before and after bifurcation (thought of as a chord) differ by one. The green and red chords have the same Maslov index, say $k$, as they are symmetric to each other, and this coincides with that of $x$ before bifurcation. This makes sure that the real Euler characteristic stays invariant. Indeed, before bifurcation we have $\mu_{L}(x)=k$ and so $\chi_L(x)=(-1)^k$; and after bifurcation, $\mu_{L}(x)=k+1$, so $\chi_L(x)=2(-1)^k+(-1)^{k+1}=2(-1)^k-(-1)^k=(-1)^k$, which explicitly shows invariance. Note that travelling through the red chord and then through the green chord gives another symmetric
periodic orbit of double period (as expected in a period-doubling bifurcation, in which the double cover $x^2$ of $x$ bifurcates). 

This can be understood within the framework provided by the GIT sequence (see Appendix \ref{app:GIT}). All monodromy matrices for different base-points along a periodic orbit are symplectically conjugated\footnote{Indeed, if $T$ is the period of the orbit $x$, and $\Psi_t$ is the Hamiltonian flow, the monodromy matrix at $x(t)$ is $D_{x(t)}\Psi_T=D_{x(0)}\Psi_{t}\circ D_{x(0)}\Psi_T\circ D_{x(t)}\Psi_{-t}=D_{x(0)}\Psi_{t}\circ D_{x(0)}\Psi_T\circ (D_{x(0)}\Psi_{t})^{-1}$.}, and hence the two reduced
monodromy matrices of $x$ at $t=0$ and $t=1/2$ induce the same element in the GIT quotient $\mathrm{Sp}(2)//\mathrm{Sp}(2)$. However, in $Sp(2)^\mathcal{I}//GL_1(\mathbb{R})$ they differ. Therefore we can apply the $B$-signature that we discussed above, to decide in practice whether the period doubling bifurcations happen at $t=0$
or $t=1/2$. Namely, the $B$-signature jumps either at $t=0$ or at $t=1/2$. The period doubling bifurcations will be symmetric at the point where 
the $B$-sign did \emph{not} jump. In particular, the $B$-signs of the negative hyperbolic critical point at the two different symmetric points of $x$ have to differ after bifurcation, precisely since only one of the points of the double cover can be symmetric, while the other one is fake symmetric. This is illustrated in Section \ref{sec:numerics} with a numerical example, in the Jupiter-Europa system.

\section{Numerics}\label{sec:numerics}

In this section, we give examples of orbit bifurcations found numerically, and illustrate the use of the various invariants discussed above. The numerical method used is the \emph{cell-mapping method}, as discussed at length in \cite{KABM}. We will focus on the (circular, restricted) three-body problem, whose Hamiltonian is given by 

$$
H: T^*\mathbb{R}^3\backslash\{M,P\}=(\mathbb R^3\setminus \{M,P\})\times \mathbb R^3\rightarrow \mathbb R,
$$
$$
H(q,p)=\frac{1}{2}\Vert p\Vert^2 - \frac{\mu}{\Vert q- M\Vert } - \frac{1-\mu}{\Vert q- P\Vert } +p_1q_2-p_2q_1, $$
where $q=(q_1,q_2,q_3)$ is the position of the Satellite, $p=(p_1,p_2,p_3)$ is its momentum, $\mu \in (0,1)$ is the mass of the secondary body M, which is fixed at $M=(\mu-1,0,0)$, and $1-\mu$ is the mass of the primary body $P=(\mu,0,0)$. The Jacobi constant is then a fixed value $c$ for $H$. The Hamiltonian $H$ is invariant under the anti-symplectic involutions 
$$
\rho: (q_1,q_2,q_3,p_1,p_2,p_3) \mapsto (q_1,-q_2,-q_3,-p_1,p_2,p_3),
$$
$$
\widetilde\rho: (q_1,q_2,q_3,p_1,p_2,p_3) \mapsto (q_1,-q_2,q_3,-p_1,p_2,-p_3),
$$
with corresponding Lagrangian fixed-point loci given by 
$$
L=\mbox{Fix}(\rho)=\{q_2=q_3=p_1=0\},
$$ 
$$
\widetilde L=\mbox{Fix}(\widetilde \rho)=\{q_2=p_1=p_3=0\}.
$$
We will then study orbits symmetric under these two symmetries, for specific cases of parameter $\mu$.

We will look at two relevant systems: the \emph{Jupiter-Europa system}, which corresponds to a circular restricted three-body problem with mass ratio $\mu=2.5266448850435e^{-05}$, and the \emph{Saturn-Enceladus system}, corresponding to $\mu=1.9002485658670e^{-07}$.

\begin{figure}
    \centering
    \includegraphics[width=0.95\linewidth]{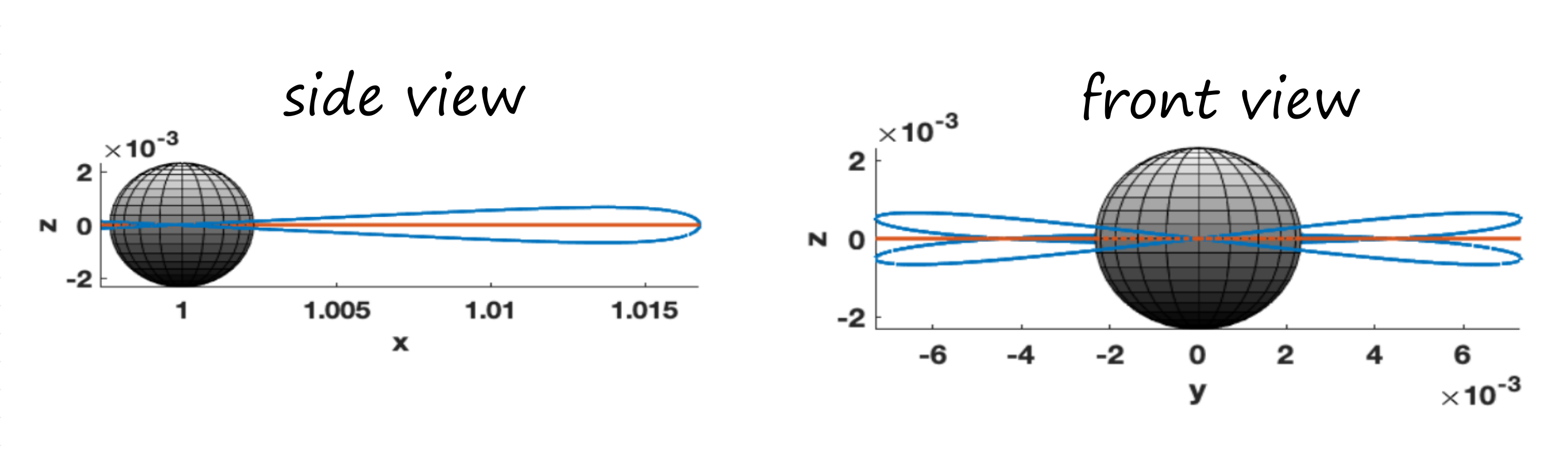}
    \includegraphics[width=0.4\linewidth]{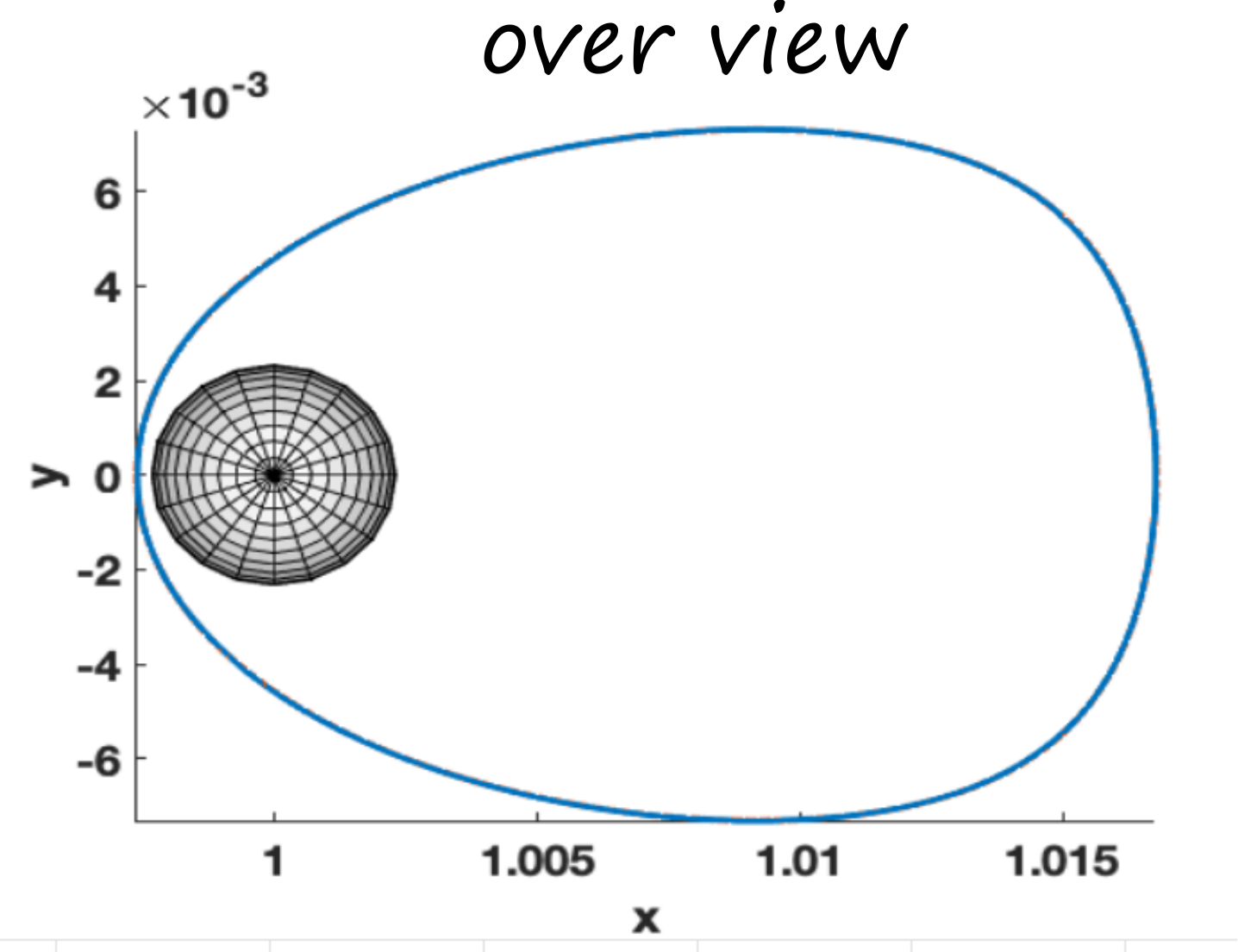}
    \includegraphics[width=0.5\linewidth]{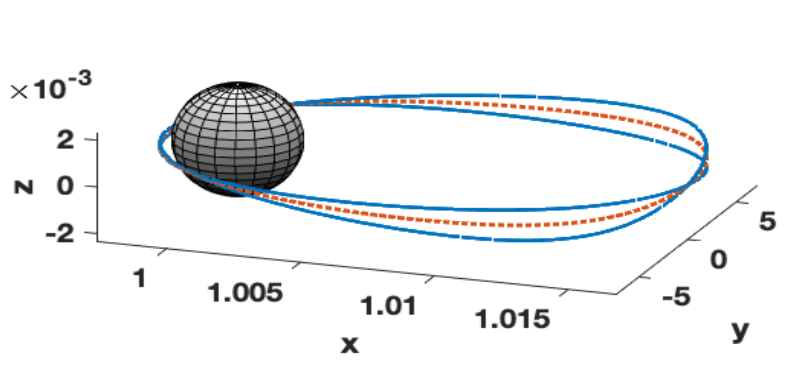}
    \caption{The prograde planar orbit $\gamma_{bef}$ (before bifurcation) is the dotted orange line; the Jacobi constant is $c=3.00357414$, and its period is $T_0=2.1215$. The spatial orbit $\beta$ of double period (after bifurcation) is the blue one; the Jacobi constant is now $c=3.003571774$, with period $T_1=4.245\approx 2T_0$ (up to small error). We call this the ``snitch'' configuration.}
    \label{fig:prograde}
\end{figure}

\begin{rem}
In what follows, the numerics were carried out in long format in MATLAB, which is higher precision as that shown here, where for readability we have truncated to 6 digits after the decimal.
\end{rem}

\begin{figure}
    \centering
    \includegraphics[width=0.6\linewidth]{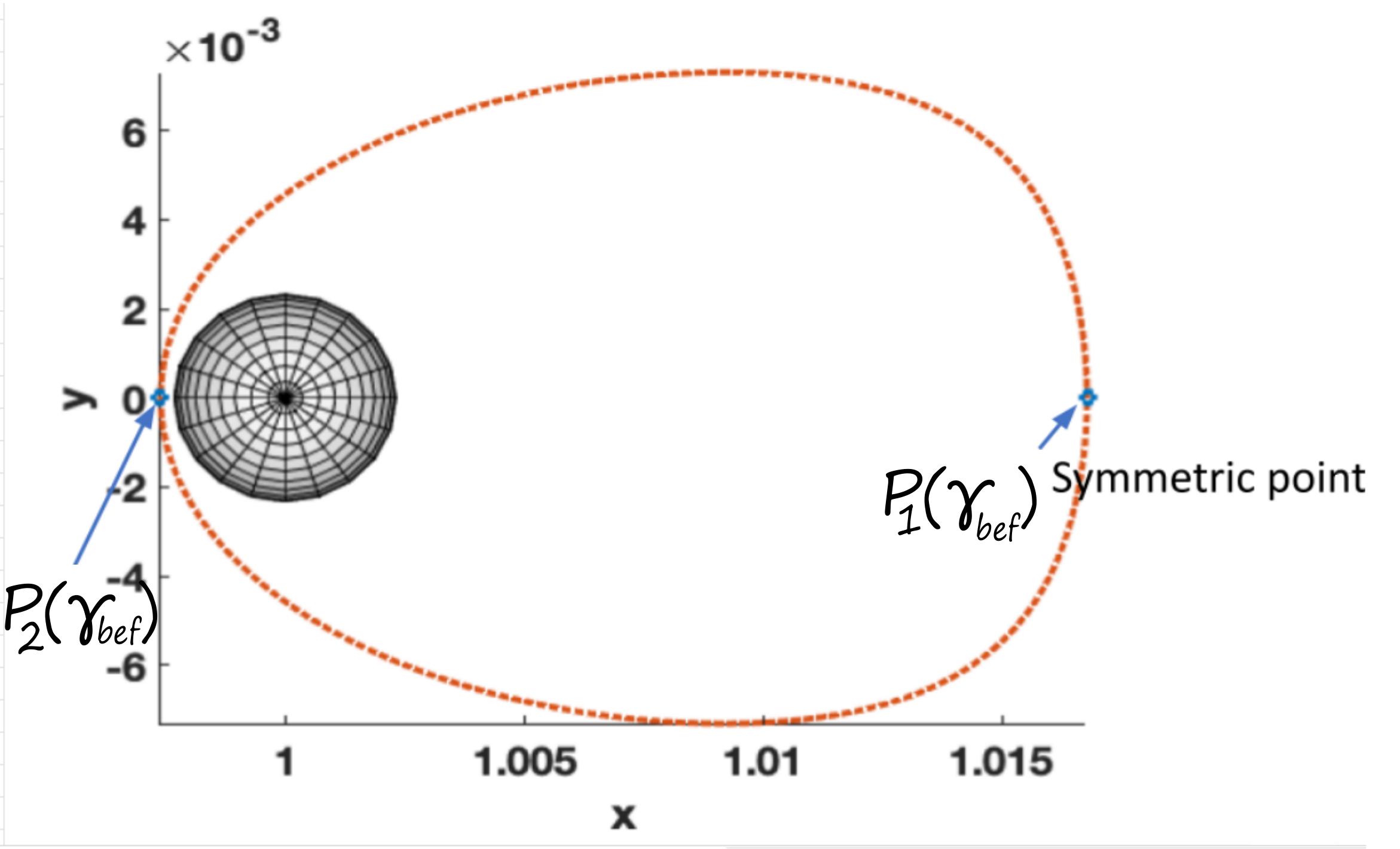}
    \caption{The two doubly symmetric points $P_1(\gamma_{bef}),P_2(\gamma_{bef})$ of the planar prograde orbit $\gamma_{bef}$.}
    \label{fig:sym_pts}
\end{figure}

\textbf{Jupiter-Europa system: period-doubling of prograde orbits ($H_2$ family)  \cite{KABM}.} As the Jacobi constant $c$ decreases, the $H_2$ family orbit depicted in Figure \ref{fig:prograde} undergoes period-doubling, i.e.\ a spatial prograde orbit of double the period appears. We denote by $\gamma_{bef}$ and $\gamma_{aft}$ the simple orbit before (i.e.\ higher energy) and after (i.e.\ lower energy) the bifurcation, and by $\beta$ the orbit with double period appearing after bifurcation. This is a \emph{doubly} symmetric period-doubling, where all orbits (i.e.\ $\gamma_{bef},\gamma_{aft}$ and $\beta$) are invariant under $\rho$ and $\widetilde \rho$. Here, $\gamma_{bef}$ is of type $\mathcal{E}^2$, and $\gamma_{aft}$ is of type $\mathcal{EH^-}$. We have, for each orbit, two $\rho$-symmetric points, where the orbit intersects $L$; similarly, two $\widetilde \rho$-symmetric ones, where the orbit intersects $\widetilde L$ (see Figure \ref{fig:sym_pts}). 

For $\gamma_{bef}$, the symmetric points are numerically found to be
$$
P_{1}(\gamma_{bef})=(1.016776,0 ,0 ,0, 0.0130372,0),
$$
$$
P_{2}(\gamma_{bef})=(0.997370,0 ,0 ,0 ,-0.125493, 0).
$$
Note that $P_1(\gamma_{bef}),P_2(\gamma_{bef})\in L\cap \widetilde L$, i.e.\ they are both $\rho$-symmetric and $\widetilde \rho$-symmetric. The (non-reduced) monodromy matrix of $\gamma_{bef}$ at $P_1(\gamma_{bef})$, is numerically computed to be:

\begin{equation*}
\centering
\begin{array}{c}\small
    M_1(\gamma_{bef})= \left(\begin{array}{ccc|ccc}
   2.930464  & 1.567115 &  0 & 0.416572 & -0.859311 & 0 \\
   -3.982191    &  -2.232667 & 0 &  -0.859311 &  1.772599 &  0\\
      0   & 0 & -0.999948 & 0 & 0 & 0.000320  \\
      \hline
        17.398921   & 6.866933 & 0 & 2.930464 & -3.982191 & 0 \\
          6.866933   & 2.056350 & 0  & 1.567115 & -2.232667 & 0 \\
            0   & 0 & -0.326763 & 0  & 0 & -0.999948 \\
\end{array}\right)\\
\end{array}
\end{equation*}

Up to small numerical rounding errors, $M_1(\gamma_{bef})$ is of the form $M_{A,B,C}$. The eigenvalues different from $1$ (which always has with multiplicity $2$), denoted $\lambda_p(\gamma_{bef})$ for \emph{planar}, and $\lambda_s(\gamma_{bef})$ for \emph{spatial}, are
$$
\lambda_{p}(\gamma_{bef})=-0.302203 + i0.953244,\; \overline{\lambda}_p(\gamma_{bef})=-0.302203 -i 0.953244,
$$
$$
\lambda_s(\gamma_{bef})=-0.999948 + i0.010225,\; \overline{\lambda}_s(\gamma_{bef})=-0.999948 -i 0.010225.
$$
Both come in elliptic conjugate pairs. Similarly, the monodromy matrix of $\gamma_{bef}$ at $P_2(\gamma_{bef})$ is 

$$
\small
M_2(\gamma_{bef})=\left(\begin{array}{ccc|ccc}
  -286.401882  & -9.995795 &  0 & 0.342004 & -9.788843 & 0 \\
   8226.036901    &  287.100358 & 0 &  -9.788864 &  280.176894 &  0\\
      0   & 0 & -0.999948 & 0 & 0 & 0.001776 \\
      \hline
       266456.859528  & 9329.998723 & 0 & -286.402124 & 8226.024432 & 0 \\
          9329.998751  & 326.685501 & 0  & -9.995804 & 287.099922 & 0 \\
            0   & 0 & -0.058878 & 0  & 0 & -0.999948 \\
\end{array}\right)
$$
Again one sees that up small errors, this is of the form $M_{A,B,C}$. By construction, $M_2(\gamma_{bef})$ is symplectically conjugated to $M_1(\gamma_{bef})$, and hence their eigenvalues need agree (we have checked that this is indeed the case, again up to small error).

\begin{figure}
    \centering
    \includegraphics[width=0.7\linewidth]{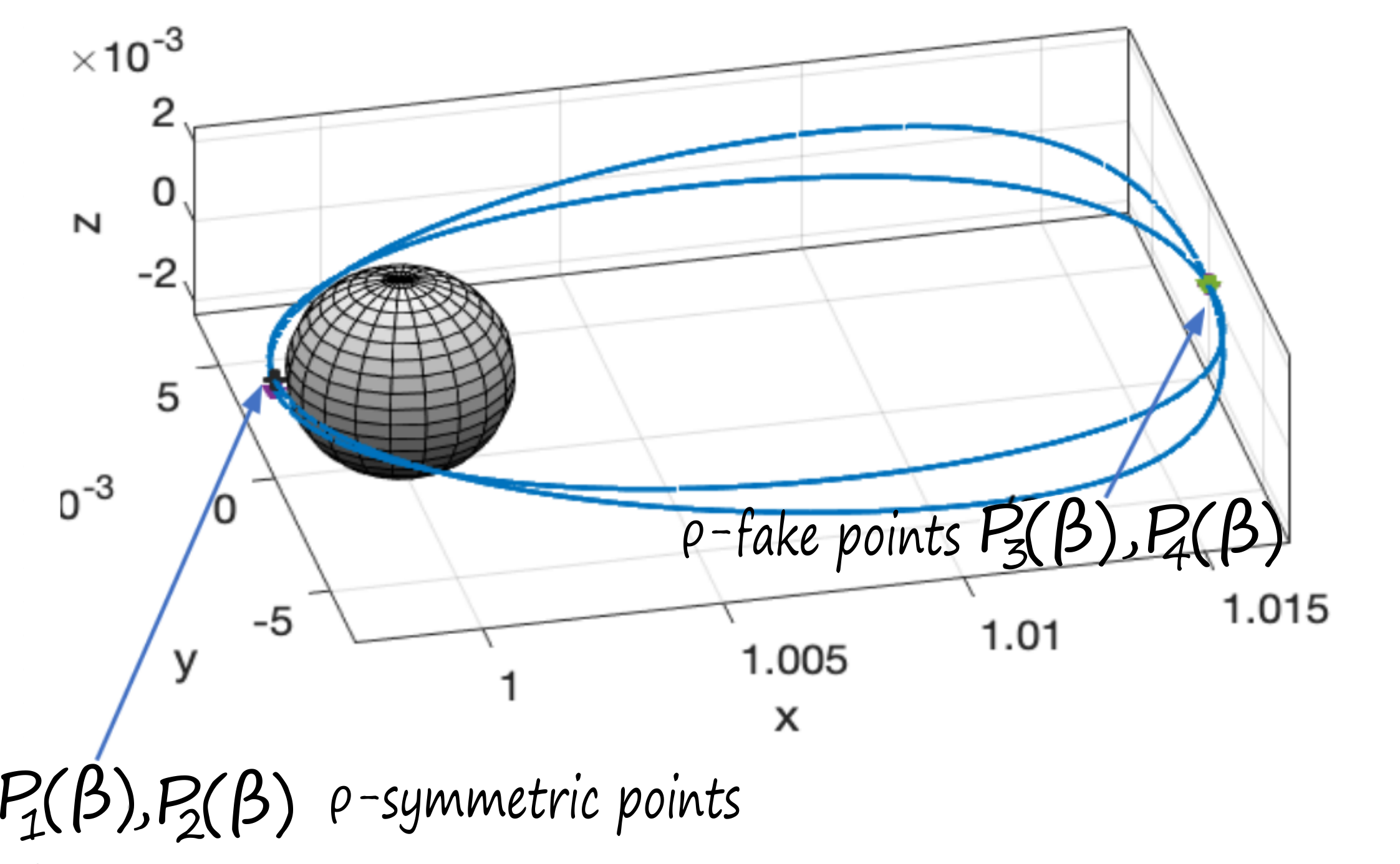}
    \caption{The $\rho$-symmetric points $P_1(\beta),P_2(\beta)$, and the $\rho$-fake points $P_3(\beta),P_4(\beta)$, of the spatial orbit $\beta$. The roles are reversed when $\rho$ is replaced by $\tilde{\rho}$.}
    \label{fig:sym_vs_fake}
\end{figure}

After bifurcation, the symmetric points of $\gamma_{aft}$ are
$$
P_{1}(\gamma_{aft})=(1.016787,0 ,0 ,0, 0.013014,0),
$$
$$
P_{2}(\gamma_{aft})=(0.997377,0 ,0 ,0 ,-0.125701, 0).
$$
Again, $P_{1}(\gamma_{aft}),P_{2}(\gamma_{aft})\in L\cap \widetilde L$ are doubly symmetric. The non-reduced monodromy matrix of $\gamma_{aft}$ at $P_1(\gamma_{aft})$ is:

\begin{equation*}
\centering
\begin{array}{c}\small
    M_1(\gamma_{aft})= \left(\begin{array}{ccc|ccc}
  2.921879 & 1.570836 &  0 & 0.412068 & -0.847786 & 0 \\
   -3.954059    &  -2.231824 & 0 &  -0.847786 &  1.744227 &  0\\
      0   & 0 & -1.000378 & 0 & 0 & -0.002449  \\
      \hline
        17.374784 & 6.880740 & 0 & 2.921879 & -3.954059 & 0 \\
          6.880740   & 2.065818 & 0  & 1.570836 & -2.231824 & 0 \\
            0   & 0 & -0.308948 & 0  & 0 & -1.000378 \\
\end{array}\right),\\
\end{array}
\end{equation*}

and that at $P_2(\gamma_{aft})$ is

\begin{equation*}
\centering
\begin{array}{c}\small
    M_2(\gamma_{aft})= \left(\begin{array}{ccc|ccc}
  -290.249559  & -10.091019 &  0 & 0.343062 & -9.857004 & 0 \\
   8368.287012    &  290.938771 & 0 &  -9.856979 &  283.215067 &  0\\
      0   & 0 & -1.000378 & 0 & 0 & 0.001672 \\
      \hline
       271816.526762   & 9480.654623 & 0 & -290.249258 & 8368.302559 & 0 \\
          9480.654594   & 330.669844 & 0  & -10.091008 & 290.939312 & 0 \\
            0   & 0 & 0.452540 & 0  & 0 & -1.000378 \\
\end{array}\right).\\
\end{array}
\end{equation*}
The eigenvalues of $M_1(\gamma_{aft})$, which up to small error coincide with that of $M_2(\gamma_{aft})$, are
$$
\lambda_{p}(\gamma_{aft})=-0.309945 + i0.950755,\; \overline{\lambda}_p(\gamma_{aft})=-0.309945 -i 0.950755,
$$
$$
\lambda_s(\gamma_{aft})=-0.972874,\; \frac{1}{{\lambda}_s(\gamma_{aft})}=-1.027883.
$$

\begin{figure}
    \centering
    \includegraphics[width=0.7\linewidth]{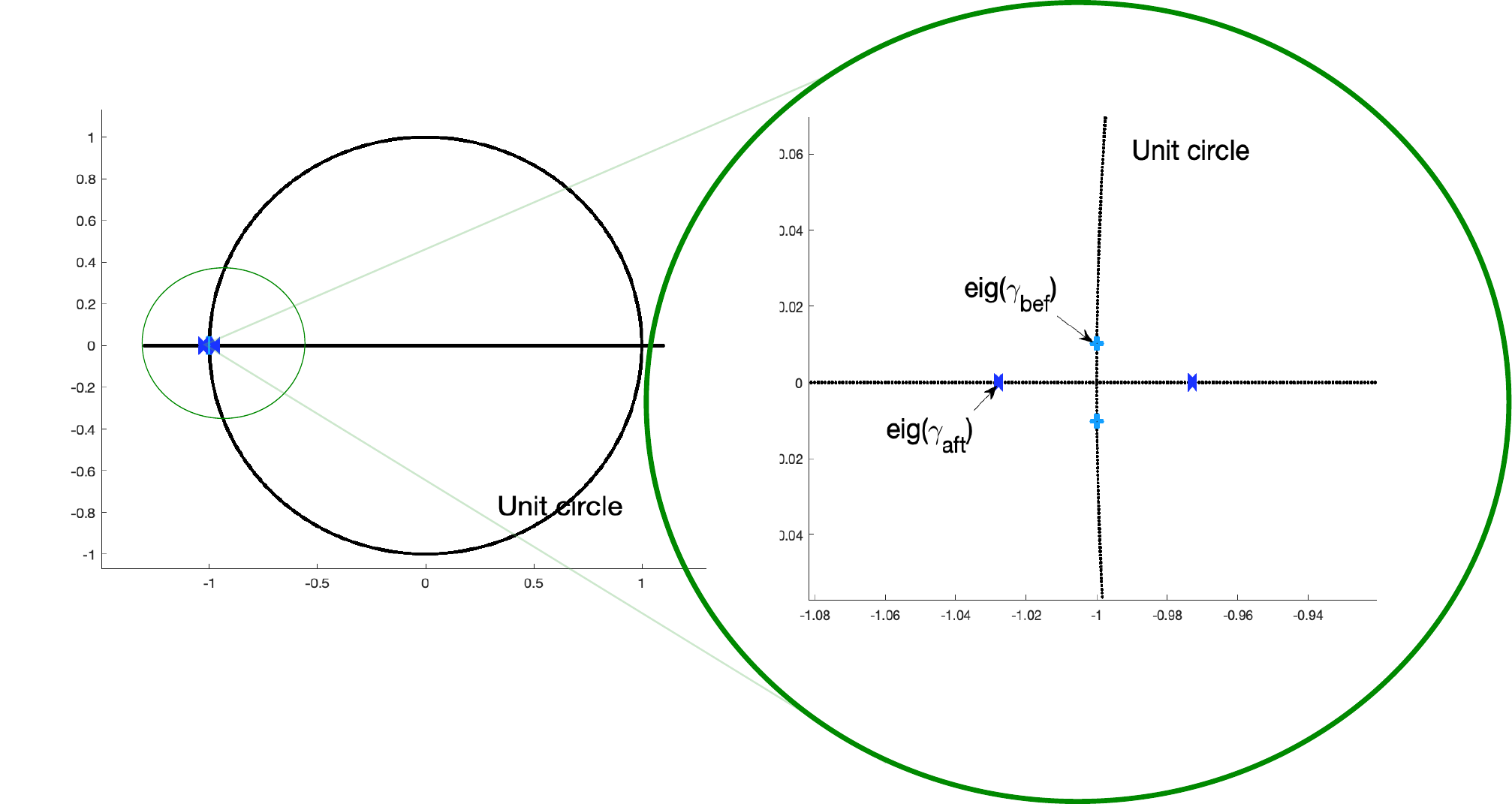}
    \caption{The spatial eigenvalues, before and after bifurcation, in two scales.}
    \label{fig:evals}
\end{figure}

Note that the planar eigenvalues stay elliptic, but the spatial ones are now a negative hyperbolic pair, as expected in a planar-to-spatial subtle division. The spatial eigenvalues, i.e.\ the bifurcating ones, are plotted in Figure \ref{fig:evals}. Now, as we discussed in the previous section, bifurcation can happen at only one of the symmetric points of $\gamma_{bef}$, for one \emph{fixed} involution. So how can we tell? Note that it is unclear just by looking at the plot in Figure \ref{fig:sym_vs_fake}. However, since in this case we have two symmetries, things become rather interesting.

We can do the same analysis for the orbit $\beta$ as we did for $\gamma_{bef}$ and $\gamma_{aft}$. Now, the algorithm used, rather than look for intersections of $\beta$ with the fixed-point loci, was implemented to look for self-intersections of $\beta$. These are: 
$$
\begin{array}{cc}
P_1(\beta)=&(0.997372,2.557023\times 10^{-17},0.000126,2.09223\times 10^{-9},-0.125462,-2.29831×\times 10^{-9}),
\\
P_2(\beta)=&(0.997372,1.768865\times 10^{-15},-0.000126,2.230434 \times 10^{-10},-0.125462,1.327282\times 10^{-8}),
\\
P_3(\beta)=&(1.016772,-3.693284\times 10^{-16},-3.907861 \times 10^{-10},2.316050 \times 10^{-9},0.013029,-0.001706),
\\
P_4(\beta)=&(1.016772,-4.489275 \times 10^{-20}, 1.590107 \times 10^{-9}, -1.587443 \times 10^{-9},0.013029,0.001706).
\end{array}
$$

But note that one cannot tell simply by inspection which of the above points are the $\rho$-symmetric ones and which ones are the $\rho$-fake ones (and similarly for $\widetilde \rho$), as they are very close to points which lie in $L$ (resp.\ $\widetilde L$), and there is numerical error involved. For this, we compute the linearizations at the corresponding points:

\begin{equation*}
\centering
\begin{array}{c}\small
    M_1(\beta)= 
\left(\begin{array}{ccc|ccc}
  -395.432864 & -13.811767 &  17.771296 & -0.211344 & 6.057227 & 0.014285 \\
   11341.304295    &  396.130733 & -508.406671 &  6.057163 &  -173.608235 &  -0.414298\\
      12.283601   & 0.426822 & 0.450370 & 0.014285 & -0.414300& -0.004380 \\
      \hline
       862435.117023  & 30173.612241 & -38842.035995 & -395.434620 & 11341.21522 & 12.283685 \\
          30173.612275   & 1055.671763 & -1358.951208 & -13.811828 & 396.127616 & 0.426825 \\
            -38842.035970   & -1358.951205 & 1749.703885 & 17.771376  & -508.402659 & 0.450367 \\
\end{array}\right)\\
\end{array}
\end{equation*}

\begin{equation*}
\centering
\begin{array}{c}\small
    M_2(\beta)= \left(\begin{array}{ccc|ccc}
  -395.431847  &-13.811731 &  -17.771251 & -0.211345 & 6.05724 & -0.014285 \\
  11341.356017    &  396.132542 & 508.409001 &  6.057139 &  -173.607555 &  0.414297 \\
      -12.283564   & -0.426821 & 0.450372 & -0.014285 & 0.414301 & -0.004380 \\
      \hline
       862435.117287  & 30173.612173 & 38842.036019 & -395.434622 & 11341.215308 & -12.283680 \\
          30173.612368   & 1055.671764 & 1358.951212  & -13.811828 & 396.127619 & -0.426825 \\
            38842.035946   & 1358.951201 & 1749.703883 & -17.771376  & 508.402662 & 0.450367 \\
\end{array}\right)\\
\end{array}
\end{equation*}

\begin{equation*}
\centering
\begin{array}{c}\small
    M_3(\beta)= \left(\begin{array}{ccc|ccc}
  3.627485  & 2.161445 &  0.112788 & -0.263329 & 0.530080 & -0.043867\\
   -5.494768    & -3.488476 & -0.147292 &  0.530087 &  -1.077505 &  0.091340 \\
      0.635462   & 0.324308 & 1.004404 & 0.043867 & -0.091341 & -0.010501 \\
      \hline
       46.989814   & 23.776561 & 0.445743 & 3.627469 & -5.494827 & -0.635463 \\
          23.776552   & 12.319463 & 0.124572  & 2.161436 & -3.488507 & -0.324309 \\
            -0.445745   & -0.124573 & 0.657275 & -0.112788  & 0.147292 & 1.004404 \\
\end{array}\right)\\
\end{array}
\end{equation*}

\begin{equation*}
\centering
\begin{array}{c}\small
    M_4(\beta)= \left(\begin{array}{ccc|ccc}
  3.6274951  & 2.161448 & -0.112789 & -0.263328 & 0.530079 & 0.043867 \\
   -5.494706   & -3.488442 & 0.147291 &  0.530093 &  -1.077514 &  -0.091339\\
      -0.635462   & -0.324308 & 1.004405 & -0.043867 & 0.091341 & -0.010501 \\
      \hline
       46.98982   & 23.776554 & -0.445741 & 3.627469 & -5.494835 & 0.635463 \\
          23.776556   & 12.31946 & -0.12457  & 2.161439 & -3.488517 & 0.324308 \\
            0.445748   & 0.124574 & 0.657275 & 0.112789  & -0.147293 & 1.004404 \\
\end{array}\right)\\
\end{array}
\end{equation*}
Now, $M_1(\beta),M_2(\beta)$ are, up to small error, of the form $M_{A,B,C}$; whereas $M_3(\beta),M_4(\beta)$ are \emph{not}. We then conclude that $P_1(\beta),P_2(\beta)$ are the $\rho$-symmetric points, whereas $P_3(\beta),P_4(\beta)$ are $\rho$-fake ones. However, these matrices implicitly assume the choice of basis, and we have chosen the basis so that $\rho$ is the standard antisymplectic involution. The roles are reversed after a change of basis for which $\tilde{\rho}$ becomes the standard such involution. After this change, one sees that $P_1(\beta),P_2(\beta)$ are $\widetilde \rho$-fake ones, and $P_3(\beta),P_4(\beta)$ are the $\widetilde \rho$-symmetric ones. So, from the perspective of $\rho$, bifurcation happened at $P_1(\beta),P_2(\beta)$, whereas from the perspective of $\tilde{\rho}$, it happened at $P_3(\beta),P_4(\beta)$. This situation is an artifact of the fact that the orbit families are doubly symmetric.

The eigenvalues of $M_1(\beta)$ (which agree with those of $M_j(\beta)$ up to small error for all $j$) are the two elliptic conjugate pairs 
$$
\lambda_{p}(\beta)=0.965396 + i0.260789,\; \overline{\lambda}_p(\beta)=0.965396 -i 0.260789,
$$
$$
\lambda_s(\beta)=-0.819634 + i0.572887,\; \overline{\lambda}_s(\beta)=-0.819634 -i 0.572887.
$$
This is of course compatible with the general discussion of symmetric period-doubling of Section \ref{sec:examples}, where the Floer numerical invariants have been used to predict what we have checked explicitly in this example. 

\smallskip

\begin{figure}
    \centering
    \includegraphics[width=0.8\linewidth]{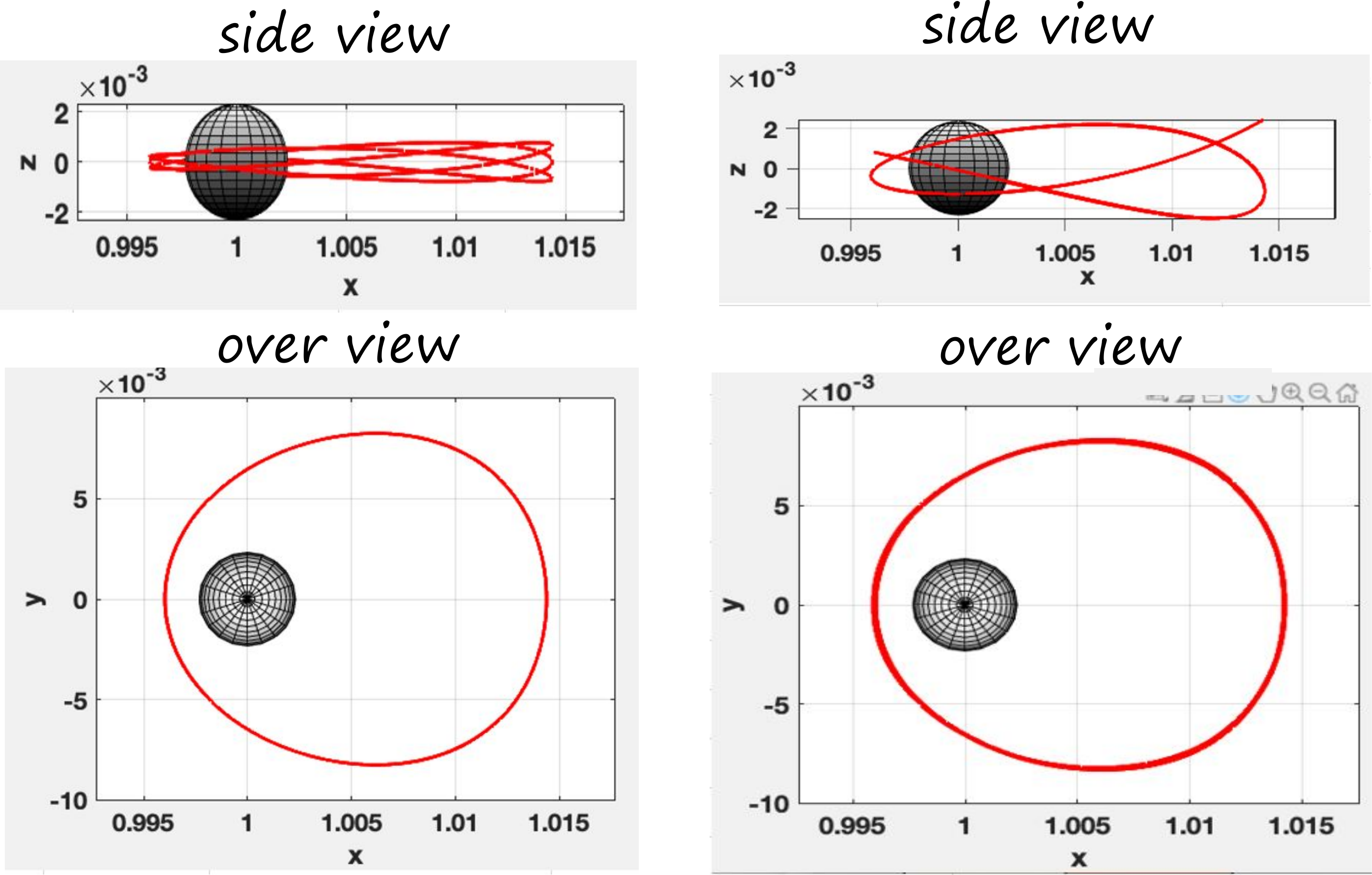}
    \caption{Left: the first period-tripling family $\gamma_1$ after bifurcation. Applying $\tilde\rho$, we obtain the family $\gamma_2$. Right: The second period-tripling family $\gamma_3$ after bifurcation. Applying $\rho$, we obtain its symmetric version $\gamma_4$.}
    \label{fig:P3bif}
\end{figure}

\textbf{Jupiter-Europa: Period tripling of prograde orbit (H$_2$ family)  \cite{KABM}.} The following is an example of a period tripling bifurcation of a $H_2$ family orbit $\gamma$, which again is doubly symmetric; see Figure \ref{fig:P3bif}. The $3$-fold cover $\gamma^3$ is of type $\mathcal{E}^2$, and bifurcates into four orbits $\gamma_1,\dots,\gamma_4$, related by the symmetries $\widetilde\rho(\gamma_1)=\gamma_2,\widetilde\rho(\gamma_3)=\gamma_4$, $\rho(\gamma_1)=\gamma_3, \rho(\gamma_2)=\gamma_4$. The orbits $\gamma_1,\gamma_2$ are in $\mathcal{E}^2$, and bifurcate from the symmetric point of $\gamma$ corresponding to $\widetilde \rho$; the orbits $\gamma_3,\gamma_4$, in $\mathcal{EH}^+$, and bifurcate from the symmetric point corresponding to $\rho$. This is compatible with the SFT-Euler characteristic: indeed, before bifurcation there is only the $3$-fold cover and so we have $\chi_{SFT}(\gamma^3)=1$. After bifurcation, the contributions of $\gamma_1,\gamma_2$ is $2$, which cancels that of $\gamma_3,\gamma_4$, which is $-2$; and we still have the contribution of the $3$-fold cover, which is $1$. So we again see that $\chi_{SFT}(\gamma^3)=1$ after bifurcation. This is also compatible with the real Euler characteristic: none of the $\gamma_i$ are symmetric, while $\gamma^3$ is; therefore $\chi_L(\gamma^3)=(-1)^{\mu_L(\gamma^3)}$ before and after.

\subsection{Numerical plots in the GIT quotient} In the following, we illustrate the numerical use of the GIT quotients via numerical plots, where we include $B$-signature computations.

\smallskip

\begin{figure}
    \centering
    \includegraphics[width=1.1\linewidth]{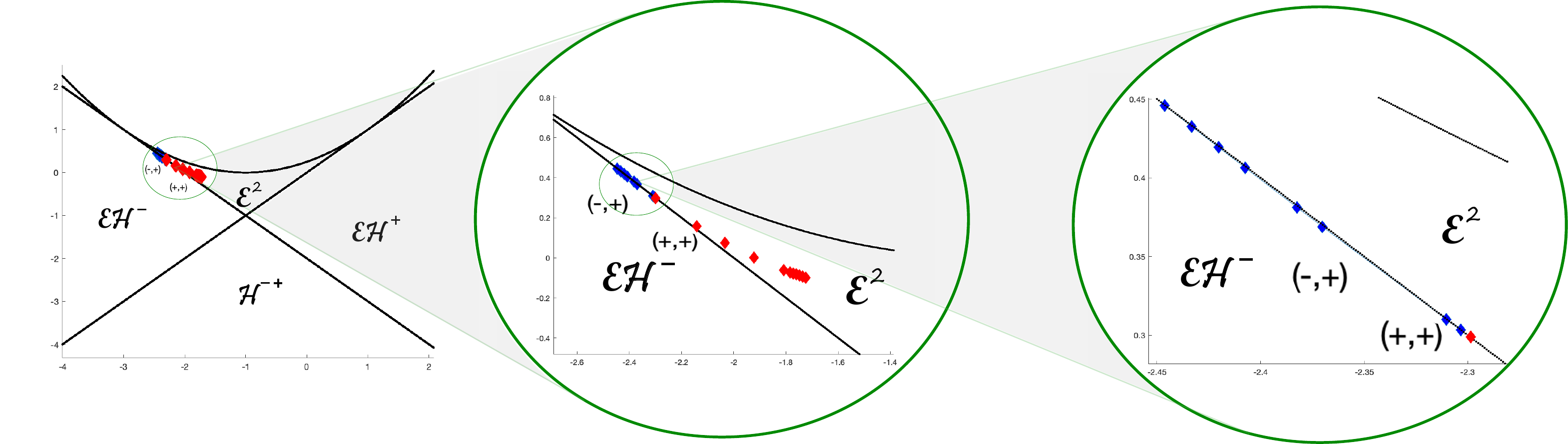}
    \caption{GIT pot of the period-doubling bifurcation of the snitch configuration.}
    \label{fig:GIT}
\end{figure}

\textbf{Snitch configuration.} We again consider the snitch configuration in the Jupiter-Europa system. Figure \ref{fig:GIT} shows a numerical plot of this period-doubling bifurcation, as seen in the base $\mathbb{R}^2$ of the GIT sequence, in three different scales. The time parameter is the Jacobi constant. Red dots correspond to $\gamma_{bef}$, and blue dots, to $\gamma_{aft}$. The bifurcation takes place when the period-doubling branch locus separating the doubly-elliptic region $\mathcal{E}^2$ and the elliptic-negative hyperbolic region $\mathcal{EH}^-$ is crossed. The plot also contains the $B$-signature of the simple orbit, before and after the bifurcation, which were computed as follows.

The $A^T$-block of $M_1(\gamma_{bef})$, at $P_1(\gamma_{bef})$, is given by
\begin{equation*}
\centering\small
    A_1^T(\gamma_{bef})= \left(\begin{array}{ccc}
   2.930464 & -3.982191 & 0 \\
   1.567115 & -2.232667 &  0\\
   0 & 0 & -0.999948
\end{array}\right)
\end{equation*}
Its two non-trivial eigenvalues are
$$
\mu_1(\gamma_{bef})=-0.999948,\;\mu_2(\gamma_{bef})=-0.302203,
$$
ordered so that $\mu_1(\gamma_{bef})<\mu_2(\gamma_{bef}),$
with corresponding eigenvectors
$$
v_1(P_1(\gamma_{bef}))=(0,0,1),v_2(P_1(\gamma_{bef}))=(0.776387,0.630256,0).
$$
Using the $B$-block of $M_1(\gamma_{bef})$, given by
\begin{equation*}
\centering\small
    B_1(\gamma_{bef})= \left(\begin{array}{ccc}
   0.416572 & -0.859311 & 0 \\
   -0.859311 &  1.772599 &  0\\
    0 & 0 & 0.00032  \\
\end{array}\right),
\end{equation*}
we simply compute
$$
\epsilon_1(P_1(\gamma_{bef}))=\mbox{sign}(v_1^T(P_1(\gamma_{bef}))\cdot B_1(\gamma_{bef})\cdot v_1( P_1(\gamma_{bef})))=\mbox{sign}(0.00032)=+,
$$
$$
\epsilon_2(P_1(\gamma_{bef}))=\mbox{sign}(v_2^T(P_1(\gamma_{bef}))\cdot B_1(\gamma_{bef}) \cdot v_2(P_1(\gamma_{bef})))=\mbox{sign}(0.114256)=+,
$$
and so the $B$-signature before bifurcation at the symmetric point $P_1(\gamma_{bef})$ is $\epsilon(P_1(\gamma_{bef}))=(+,+)$ (as depicted in Figure \ref{fig:GIT}). The same procedure applied to $M_1(\gamma_{aft})$ gives
$$
\mu_1(P_1(\gamma_{aft}))=-1.00038,\; \mu_2(P_1(\gamma_{aft}))=-0.309942,  
$$
with corresponding eigenvectors
$$
v_1(P_1(\gamma_{aft}))=(0,0,1),\;v_2(P_1(\gamma_{aft}))=(0.774275,0.632849,0),
$$
and corresponding $B$-signs
$$
\epsilon_1(P_1(\gamma_{aft}))=\mbox{sign}(v_1^T(P_1(\gamma_{aft}))\cdot B_1(\gamma_{aft})\cdot v_1(P_1(\gamma_{aft})))=\mbox{sign}(-0.00245)=-,
$$
$$
\epsilon_2(P_1(\gamma_{aft}))=\mbox{sign}(v_2^T(P_1(\gamma_{aft}))\cdot B_1(\gamma_{aft}) \cdot v_2(P_1(\gamma_{aft})))=\mbox{sign}(0.114766)=+,
$$
and so the $B$-signature after bifurcation at the symmetric point $P_1(\gamma_{aft})$ is $\epsilon(P_1(\gamma_{aft}))=(-,+)$ (also depicted in Figure \ref{fig:GIT}). 

We check explicitly in this example the fact, alluded to in the general discussion of Section \ref{sec:examples}, that the $B$-signature at different symmetric points of will differ after bifurcation in a symmetric period-doubling. Indeed, replacing $P_1$ by $P_2$, the eigenvalues of the corresponding $A^T$-blocks $A_2^T(\gamma_{bef}),A_2^T(\gamma_{aft})$ need respectively coincide with the $\mu_i(\gamma_{bef}),\mu_i(\gamma_{aft})$ for $i=1,2$ (checked up to  numerical error), and the corresponding eigenvectors are  
$$
v_1(P_2(\gamma_{bef}))=(0,0,1),\; v_2(P_2(\gamma_{bef}))=(-0.999396,-0.0347591,0). 
$$
$$
v_1(P_2(\gamma_{aft}))=(0,0,1),\; v_2(P_2(\gamma_{aft}))=(-0.9994,-0.0346264,0),
$$
with associated $B$-signs
$$
\epsilon_1(P_2(\gamma_{bef}))=\mbox{sign}(v_1^T(P_2(\gamma_{bef}))\cdot B_2(\gamma_{bef})\cdot v_1(P_2(\gamma_{bef})))=\mbox{sign}(0.001776 )=+,
$$
$$
\epsilon_2(P_2(\gamma_{bef}))=\mbox{sign}(v_2^T(P_2(\gamma_{bef}))\cdot B_2(\gamma_{bef}) \cdot v_2(P_2(\gamma_{bef})))=\mbox{sign}(6.86473\times 10^{-6})=+,
$$
$$
\epsilon_1(P_2(\gamma_{aft}))=\mbox{sign}(v_1^T(P_2(\gamma_{aft}))\cdot B_2(\gamma_{aft})\cdot v_1(P_2(\gamma_{aft})))=\mbox{sign}(0.001672 )=+,
$$
$$
\epsilon_2(P_2(\gamma_{aft}))=\mbox{sign}(v_2^T(P_2(\gamma_{aft}))\cdot B_2(\gamma_{aft}) \cdot v_2(P_2(\gamma_{aft})))=\mbox{sign}(7.10883\times 10^{-6})=+.
$$
We see that the $B$-signatures are $\epsilon(P_2(\gamma_{bef}))=\epsilon(P_2(\gamma_{aft}))=(+,+)$, and therefore the $B$-sign of the eigenvalue $\mu_1$ (the one undergoing bifurcation) indeed differs after bifurcation, for different choices of symmetric point. As pointed out in the general case, the fact that there is a sign jump at $P_1$ and not at $P_2$ indicates that $P_2$ gives rise to the $\rho$-symmetric points, and not $P_1$, where $\rho$ is the involution which is standard in the current choice of basis.

We conclude this section with a series of plots, including examples in the Saturn-Enceladus system, some of which are also discussed in \cite{KABM} and \cite{KABM2018}.

\begin{figure}
    \centering
    \includegraphics[width=0.9\linewidth]{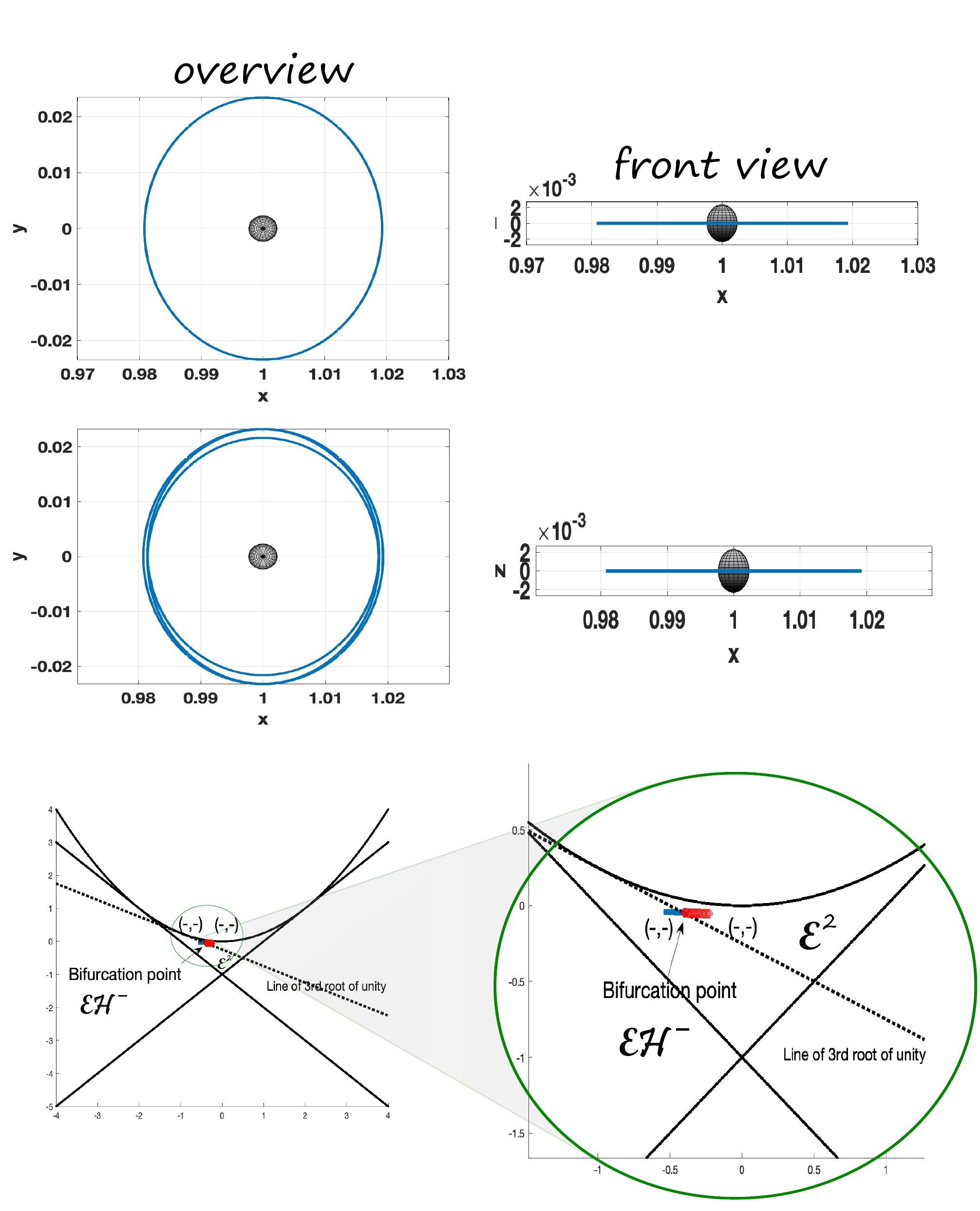}
    \caption{Jupiter-Europa system: a symmetric planar to planar period-tripling bifurcation of DRO family. Above: the planar simple orbit at bifurcation, the \emph{distant retrograde orbit}, $c=2.9999$, $T_0=2.504$. Middle: the planar triple period orbit after bifurcation $c\gtrsim 2.9999$, $T= 7.3\approx 3T_0$. Below: GIT plot, including $B$-signs.}
    \label{fig:JEDRO3}
\end{figure}

\begin{figure}
    \centering
    \includegraphics[width=0.9\linewidth]{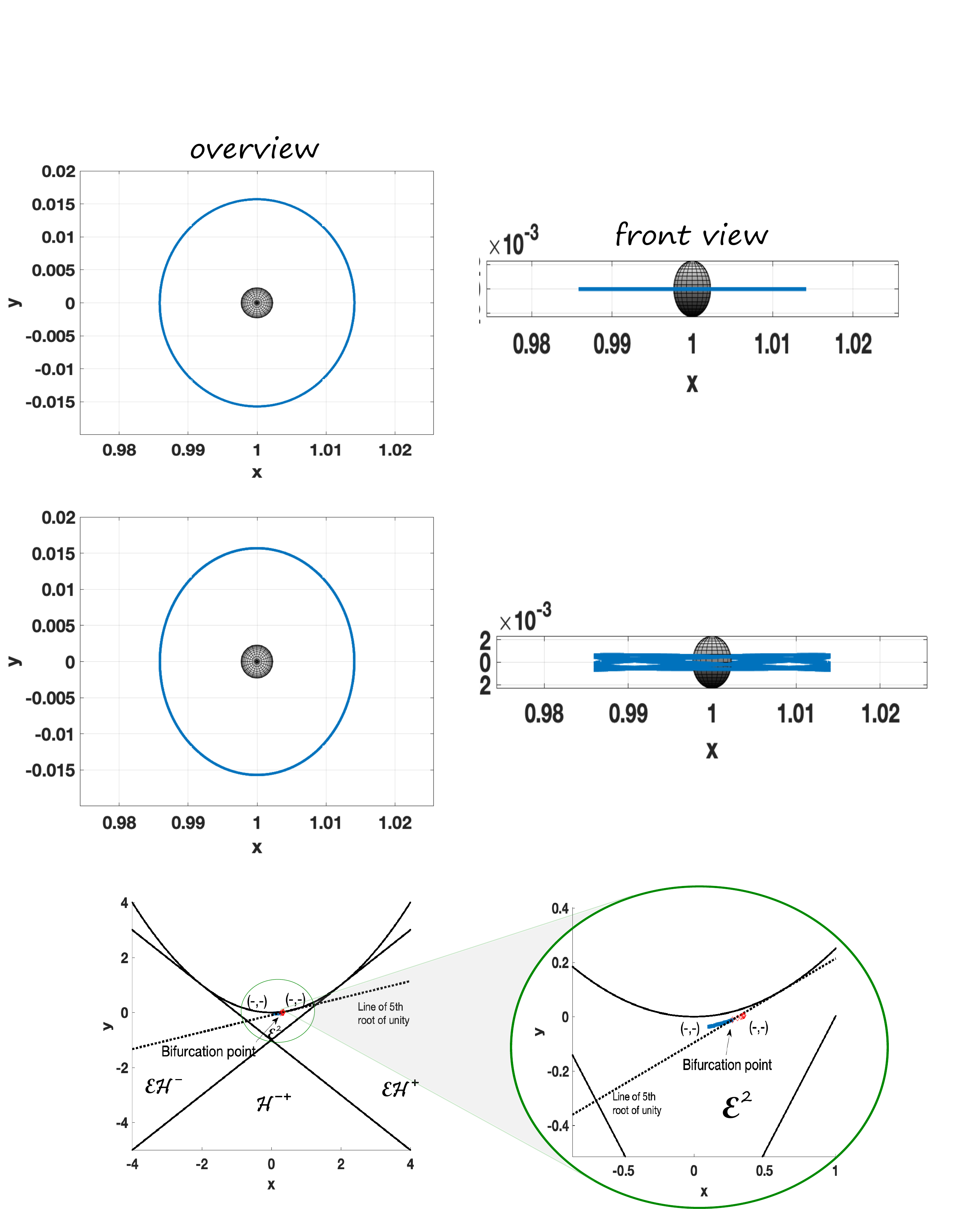}
    \caption{Jupiter-Europa system: a symmetric planar to spatial $5$-fold bifurcation of DRO family. Above: the planar simple orbit at bifurcation, a distant retrograde orbit in the same family of Figure \ref{fig:JEDRO3}, but with  $c=3.0005,$ $T_0=1.705$. Middle: the spatial $5$-fold period orbit after bifurcation, $c\gtrsim 3.0005,$ $T=8.52\approx 5T_0$. Below: GIT plot, including $B$-signs.}
\end{figure}

\begin{figure}
    \centering
    \includegraphics[width=0.9\linewidth]{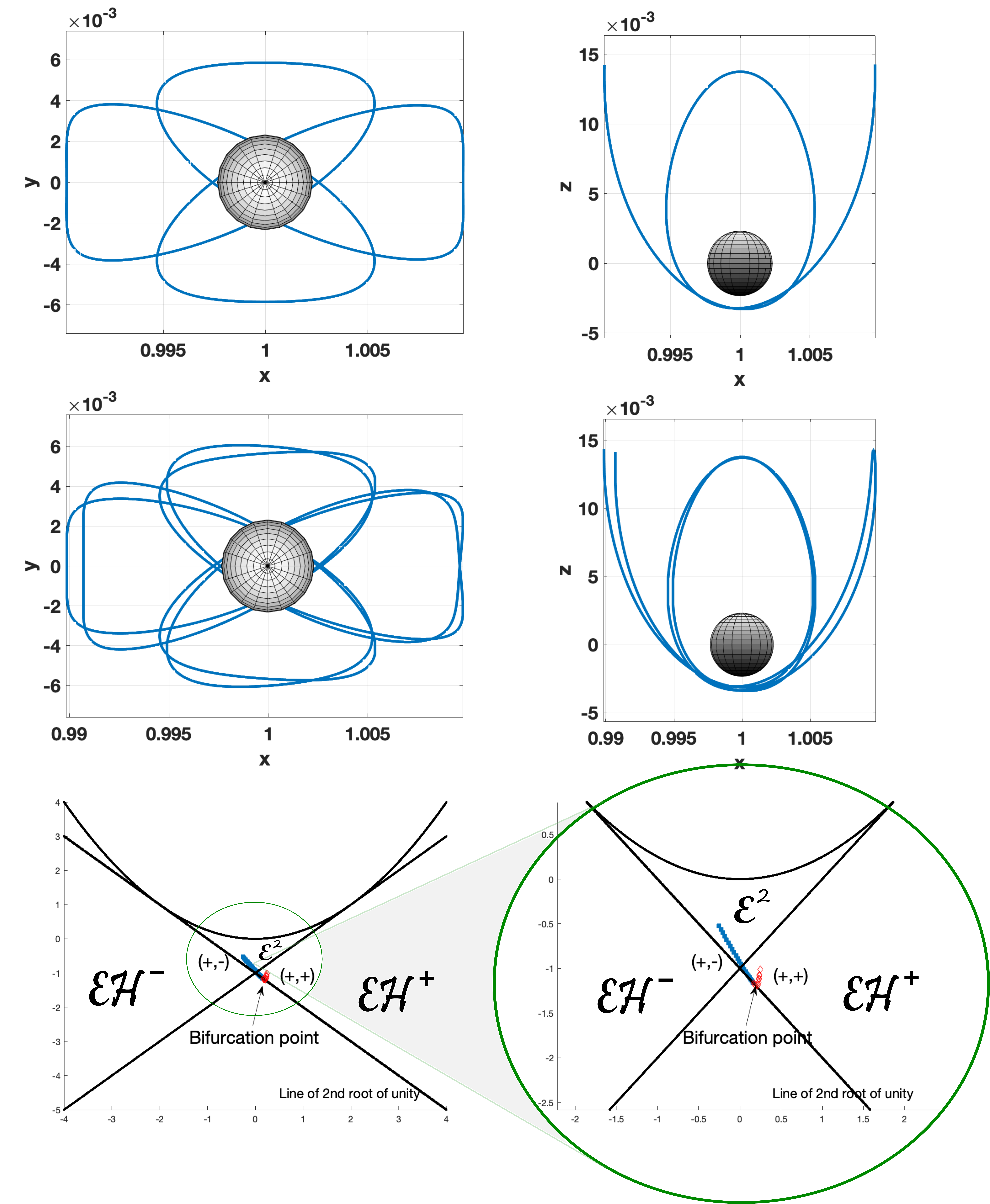}
    \caption{Jupiter-Europa system: a symmetric spatial to spatial period-doubling bifurcation. Above: the spatial simple orbit at bifurcation, $c=3.0028$, $T_0=4.62$. Middle: the spatial period-doubling orbit after bifurcation $c\gtrsim 3.0028,$ $T=9.23\approx 2T_0$. Below: GIT plot, including $B$-signs.}
\end{figure}

\begin{figure}
    \centering
    \includegraphics[width=0.86\linewidth]{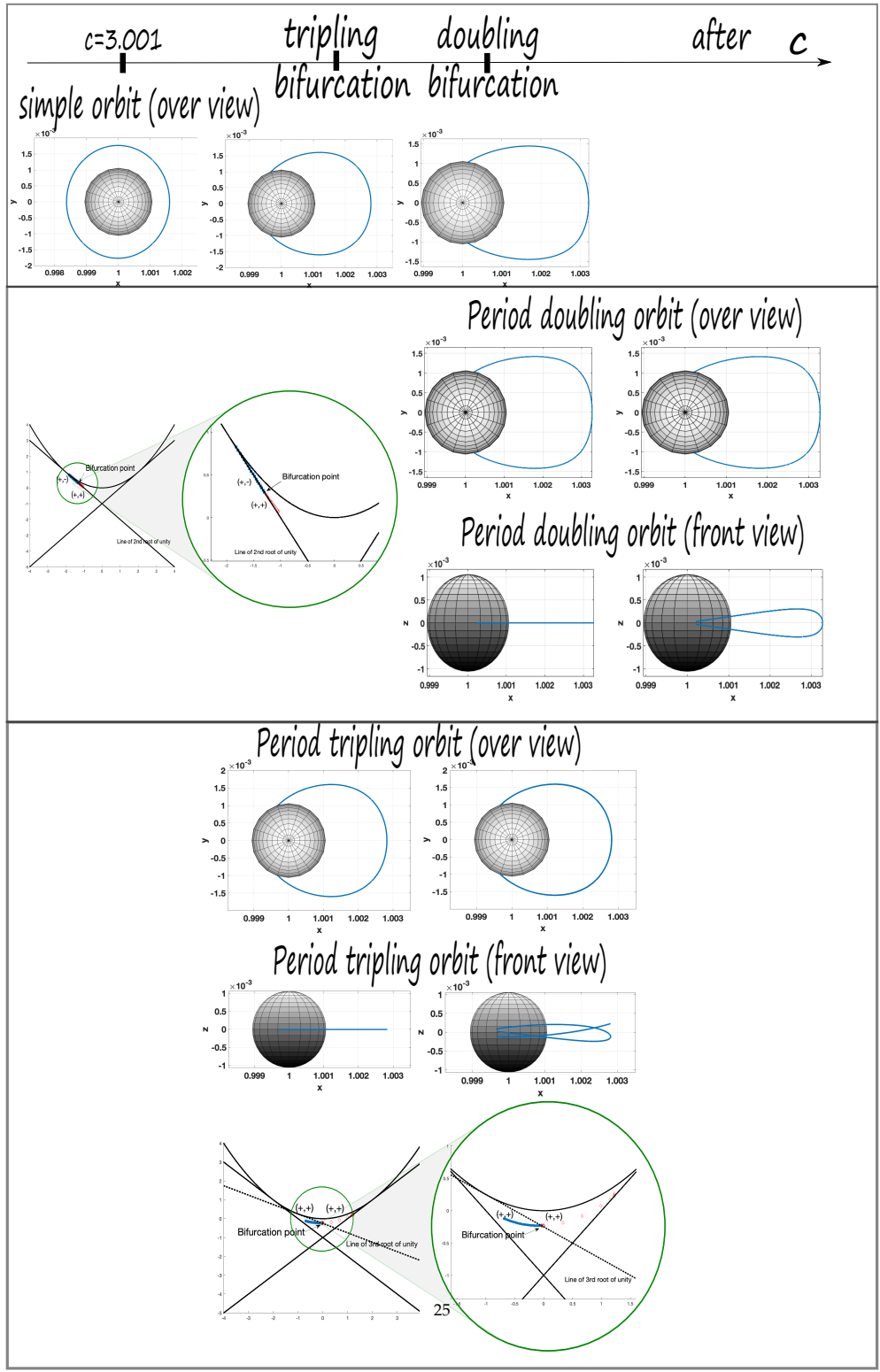}
    \caption{Saturn-Enceladus system: two symmetric planar to spatial bifurcations of the same family of planar orbits, one period-doubling, and one period-tripling. Above: $T_1=1.2, T_2=1.6, T_3=2$ respectively. Middle: $T=4.2\approx 2T_3$ after bifurcation. Below: $T=4.85\approx 3T_2$ after bifurcation.}
    \label{fig:my_label}
\end{figure}

\clearpage

\appendix

\section{GIT sequence}\label{app:GIT}

\begin{figure}
    \centering
    \includegraphics[width=0.55 \linewidth]{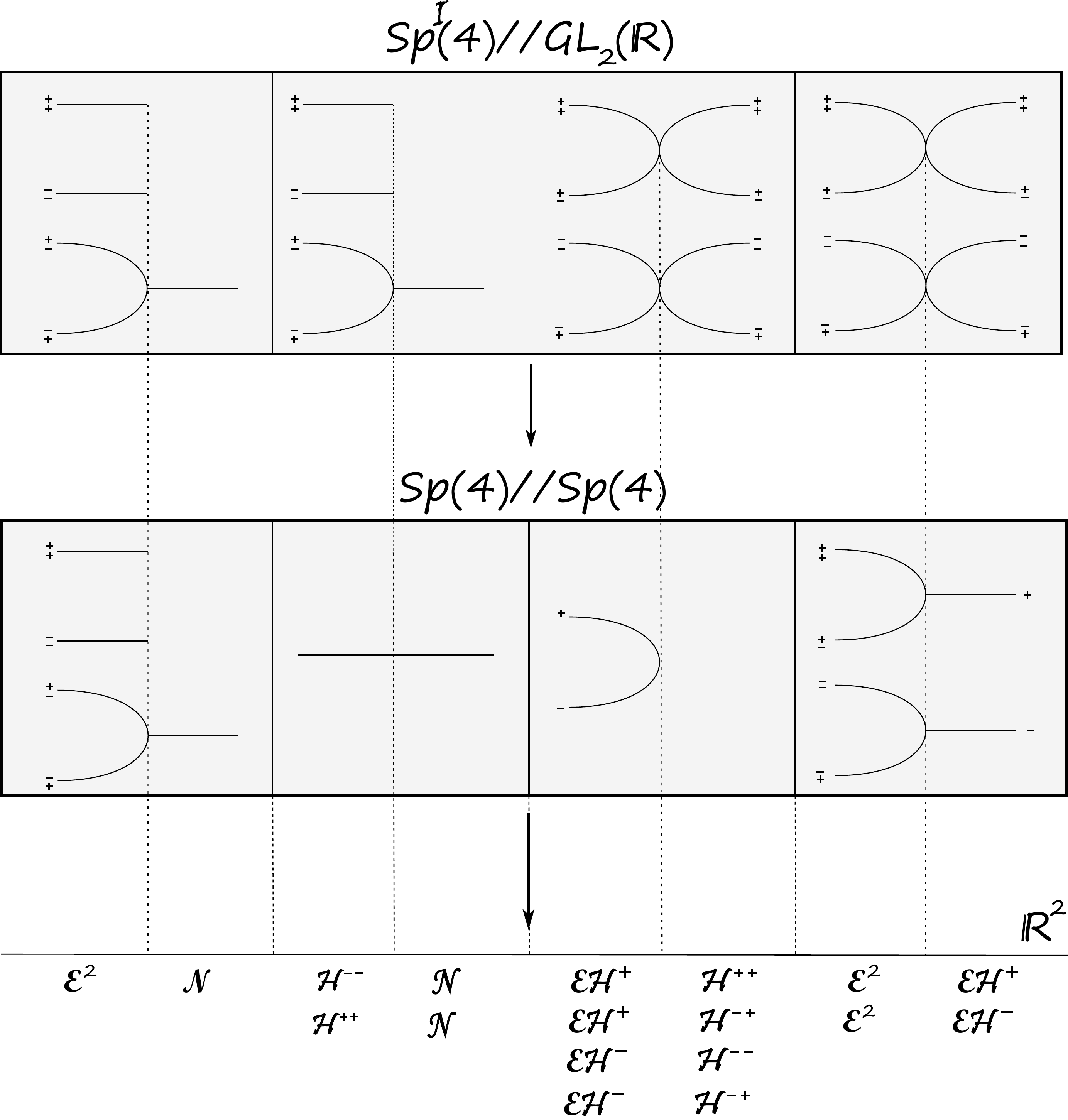}
    \caption{This picture shows the branches of $\mathrm{Sp}^\mathcal{I}(4)//\mathrm{GL}_2(\mathbb{R})$ and $\mathrm{Sp}(4)//\mathrm{Sp}(4)$, which are $2$-dimensional ``sheets'' covering the different regions of the plane depicted in Figure \ref{fig:basen2} (we drop one dimension for visualization). The signs on each branch correspond to $B$-positivity/negativity of the corresponding eigenvalues (a priori there are $4$ possibilities, since there are two eigenvalues). The first vignette shows how they come together when crossing from $\mathcal{E}^2$ to $\mathcal{N}$. On the second, when crossing from $\mathcal{H}^{--}$ to $\mathcal{N}$; the picture is the same for $\mathcal{H}^{++}$ to $\mathcal{N}$, and so on. All branches come together to a single point along each of the three singular points $(2,1),(0,-1), (-2,1)$. The map $\mathrm{Sp}^\mathcal{I}(4)//\mathrm{GL}_2(\mathbb{R})\rightarrow \mathrm{Sp}(4)//\mathrm{Sp}(4)$ in the GIT sequence collapses branches together, as shown. For example, $B$-positivity/negativity over the hyperbolic eigenspace of matrices of type $\mathcal{EH}^+$ is not invariant under symplectic conjugation, and hence the corresponding branches come together in $\mathrm{Sp}(4)//\mathrm{Sp}(4)$.}
    \label{fig:bifurc}
\end{figure}

We now explain the notion of the GIT sequence. The rough idea is to understand the topology of the configuration space consisting of the collection of all possible pairs $(p,\epsilon)$ of point $p=(\det(A),\mbox{tr}(A))$ and $B$-signature $\epsilon=(\epsilon_1,\epsilon_2)$, together with the structure of the projection $(p,\epsilon)\mapsto p \in \mathbb{R}^2$. This is illustrated in Figure \ref{fig:bifurc}, where the configuration space lies on the top, and has different ``branches'' corresponding to different $B$-signature $\epsilon$, which get collapsed on top of each other under the projection. The plus/minus labels in the branches of the ``middle'' space of Figure \ref{fig:bifurc} records precisely the Krein sign over elliptic components. Whenever an orbit is symmetric and we choose a symmetric point, it may be ''lifted'' uniquely from the ''middle'' space to the top one. One can do this more formally, as follows. The treatment will assume some mathematical background, and is included for completeness. 

\begin{rem}[GIT quotient] To give the definition of the GIT sequence, we need to introduce some terminology. Recall that if a group $G$ acts on a topological space $X$, the geometric quotient $X/G$ is the space of $G$-orbits, i.e.\ a point in $X/G$ is a set of the form $\{g\cdot x:x\in X\}\subset X$. In general, this space might not be Hausdorff, i.e.\ there might be points which cannot be separated from other points. To fix this, one considers the \emph{GIT quotient}, the space $X//G$ obtained by identifying two points of $X$ if the closures of their $G$-orbits intersect (and this space is indeed Hausdorff). We shall consider only GIT quotients in what follows, although the reader might choose to ignore this technicality. 
\end{rem}

The GIT sequence consists of the sequence of maps \begin{equation}\label{seq}
\mathrm{Sp}^\mathcal{I}(2n)//\mathrm{GL}_n(\mathbb{R}) \rightarrow \mathrm{Sp}(2n)//\mathrm{Sp}(2n)
\rightarrow \mathrm{M}_{n\times n}(\mathbb{R})//\mathrm{GL}_n(\mathbb{R})\cong \mathbb R^n,
\end{equation}
given by
$$[M_{A,B,C}] \mapsto [[M_{A,B,C}]] \mapsto [A].$$ Here,
$\mathrm{Sp}(2n)$ is the symplectic group, which acts on itself by conjugation, i.e.\ via $A\cdot B=A B A^{-1}$; and $\mathrm{GL}_n(\mathbb{R})$ also acts by conjugation on the space of matrices $\mathrm{M}_{n\times n}(\mathbb{R})$. Above we denote by $[M_{A,B,C}]$ the equivalence class of the matrix
$M_{A,B,C} \in \mathrm{Sp}^\mathcal{I}(2n)$ in the GIT quotient $\mathrm{Sp}^\mathcal{I}(2n)//\mathrm{GL}_n(\mathbb{R})$, by $[[M_{A,B,C}]]$ the equivalence class in the GIT quotient $\mathrm{Sp}(2n)//\mathrm{Sp}(2n)$, and by $[A]$ the equivalence class of the first block $A\in M_{n\times n}(\mathbb{R})$ in $ \mathrm{M}_{n\times n}(\mathbb{R})//\mathrm{GL}_n(\mathbb{R})$. We have used the fact that mapping the equivalence class of a matrix $A \in \mathrm{M}_{n\times n}
(\mathbb{R})$ to the coefficients of its characteristic polynomial, we get an identification $\mathrm{M}_{n\times n}(\mathbb{R})//\mathrm{GL}_n(\mathbb{R}) \cong \mathbb{R}^n$; see \cite[Appendix A]{FM}.

In the examples above, where the spaces consists of matrices, the transition from to the geometric quotient to the GIT quotient basically means, in practice, to ignore Jordan factors, replacing them with diagonal blocks. The resulting matrices, while not
necessarily equivalent in the original quotient, become so in the GIT one, see \cite[Appendix A]{FM}. In \cite{FM}, the cases $n=1$ and $n=2$ (for instance, relevant for the planar and the spatial three-body problems, respectively) are studied in detail. In particular, the topology of these GIT quotients is fully determined, as well as the maps. In this article, we will make use of the case $n=2$, where the base of the GIT sequence is the plane $\mathbb{R}^2$, together with the structure of its bifurcation loci, as shown in Figure \ref{fig:complete_bifurcation}. The maps of the sequence are also very concrete and therefore simple to implement, i.e.\ given by
$$
\mathrm{Sp}^{\mathcal{I}}(4)//\mathrm{GL}_2(\mathbb{R})\rightarrow Sp(4)//Sp(4) \rightarrow \mathbb{R}^2
$$
$$
[M_{A,B,C}]\mapsto [[M_{A,B,C}]]\mapsto (\mathrm{tr}(A),\det(A))=p,
$$
which motivates the construction of $p$ that we explained above. Indeed, as follows from \cite{FM}, the GIT quotient $\mathrm{Sp}^{\mathcal{I}}(4)//\mathrm{GL}_2(\mathbb{R})$ is precisely the configuration space for the pairs $(p,\epsilon)$.

\section{Invariance of the SFT-Euler characteristic}\label{app:invariance}

While in general invariance of $\chi_{SFT}$ is non-trivial, for the \emph{generic} bifurcations in dimension four, we can check it directly, as follows.

\subsection{Generic four-dimensional bifurcations} We follow the cases as listed in the book by Abraham and Marsden \cite{abraham-marsden}.

\smallskip

\textbf{Creation \cite[p.\,598]{abraham-marsden}: } In this case initially there was no periodic orbit
at all. Hence $\chi_{SFT}=0.$ After the creation there is a simple elliptic and positive hyperbolic orbit. In particular, the SFT-Euler characteristic stays zero. 
\\ \\
\textbf{Subtle division \cite[p.\,599]{abraham-marsden}: }In this case the double cover of an elliptic orbit bifurcates. We consider the SFT-Euler characteristic for the simple
orbit, which we denote by $\chi_{SFT}^1$, as well as for the double cover, denoted $\chi_{SFT}^2$. Before
the transition there is one simple elliptic orbit. Therefore $\chi_{SFT}^1=-1.$ After the transition the simple orbit becomes negative hyperbolic. There is no bifurcation of the simple
periodic orbit, just its double cover bifurcates. Hence $\chi_{SFT}^1$ stays minus one. For invariance of $\chi_{SFT}^2$, note that the double cover of an elliptic orbit is elliptic as well. Therefore $\chi^2_{SFT}=-1.$ After the transition the simple elliptic orbit becomes negative hyperbolic. Its double cover is therefore a bad positive hyperbolic orbit and does not contribute to the SFT-Euler characteristic. The orbit which bifurcates is elliptic and hence the SFT-Euler characteristic stays $-1$. 
\\ \\
\textbf{Murder \cite[p.\,600]{abraham-marsden}: } We consider the SFT-Euler characteristic $\chi_{SFT}^1$ of the simple orbit as well as the one for the double
cover $\chi_{SFT}^2$. The case for the simple orbit is completely analogous as in the subtle division.
An elliptic periodic orbit becomes negative hyperbolic and therefore $\chi_{SFT}^1=-1.$ However, the case of the double cover is different. Here before bifurcation we have a double covered elliptic one and a simple positive hyperbolic one. A simple positive hyperbolic orbit is good and the double cover of an elliptic orbit is elliptic as well. Therefore $\chi_{SFT}^2=0.$
After bifurcation just the double cover of the negative hyperbolic orbit is left. This is a bad positive
hyperbolic orbit and therefore does not contribute to the SFT-Euler characteristic. Then $\chi_{SFT}^2$ stays zero after the transition. 
\\ \\
\textbf{Phantom kiss \cite[p.\,602]{abraham-marsden}: } We discuss the 3-kiss, via the SFT-Euler characteristic of the 3-fold cover $\chi_{SFT}^3$. Before the bifurcation we have a 3-fold covered elliptic orbit and a simple positive hyperbolic one. Therefore $\chi_{SFT}^3=0.$ After bifurcation we still have a 3-fold covered elliptic orbit and a positive hyperbolic one, so that the SFT-Euler characteristic does not change. The discussion for the 4-kiss
is similar, when one considers the SFT-Euler characteristic of the 4-fold cover $\chi_{SFT}^4$. 
\\ \\
\textbf{Emission \cite[p.\,603]{abraham-marsden}: } We discuss here the case $p=4$ as illustrated in the figure in \cite[p.\,603]{abraham-marsden}. We consider the SFT-Euler characteristic for the 4-fold cover
$\chi_{SFT}^4$. Before bifurcation there is one 4-fold covered elliptic orbit. Therefore $\chi_{SFT}^4=-1.$
After bifurcation there is a 4-fold covered elliptic orbit, a simple elliptic orbit and a simple
positive hyperbolic orbit. We see again that the SFT-Euler characteristic does not change. 

\subsection{Non-generic four-dimensional bifurcations.} There are relevant problems in celestial mechanics where there are bifurcations which do not fall in the generic classification. We now discuss some of them.

\smallskip

\textbf{Hénon families in Hill's lunar problem.} Although in theory the probability to have a non-generic bifurcation is basically zero, in practice non-generic bifurcations occur quite often. The reason is that the Hamiltonians one usually considers are invariant under various symmetries. A nongeneric bifurcation was described by H\'enon 
in \cite{henon} while studying Hill's lunar
problem, a limit case of the restricted three-body problem where the massless body is assumed very close to the small primary. Hill's lunar problem can therefore be considered as an approximation to
the Jupiter-Europa or Saturn-Enceladus systems, when one lets the mass of Europa, respectively Enceladus, go to zero. While the potential of the restricted three-body problem is invariant
under reflection at the $x$-axis, i.e.\ the axis on which the two primaries
lie, Hill's lunar problem is additionally invariant under reflection at the
$y$-axis. The family of the direct or prograde periodic orbit is referred to as family $g$. The direct orbit
is invariant under reflection at the $x$-axis as well as under reflection at the $y$-axis. For small energy the direct orbit is elliptic. However, for higher energy it becomes positive hyperbolic. At the bifurcation point two new families, referred to as $g'$, appear. These two families
are still invariant under reflection at the $x$-axis but not anymore under reflection at the $y$-axis. Instead of that, reflection at the
$y$-axis maps one branch of the $g'$-family to the other branch. As explained by H\'enon \cite{henon}, at their birth, the two $g'$-branches
are elliptic. 

We can now check the invariance of the SFT-Euler characteristic
for this nongeneric bifurcation. Before the bifurcation the direct orbit
was elliptic. Therefore the SFT-Euler characteristic is minus one. After
the bifurcation the direct periodic orbit is positively hyperbolic. Since
it is simple it is a good positive hyperbolic orbit and therefore contributes
$+1$ to the SFT-Euler characteristic. However, after bifurcation we have
to take into account in addition the two $g'$-periodic orbits which
are both elliptic and therefore contribute each $-1$ the the SFT-Euler characteristic. So their sum $1-1-1=-1$ stays minus one.


\begin{thebibliography}{99}
 \bibitem{abraham-marsden}R.\,Abraham, J.\,Marsden,
 \emph{Foundations of Mechanics}, 2nd ed. Addison-Wesley, New York
 (1978).
 
 \bibitem{Broucke} R.\ Broucke, Stability of periodic orbits in the elliptic, restricted three-body problem. AIAA J. 7,1003 (1969).
 
 \bibitem{EGH} Eliashberg, Y.; Givental, A.; Hofer, H. Introduction to symplectic field theory. GAFA 2000 (Tel Aviv, 1999). Geom. Funct. Anal. 2000, Special Volume, Part II, 560--673.
 
 \bibitem{F1} Floer, Andreas. A relative Morse index for the symplectic action. Comm. Pure Appl. Math. 41 (1988), no. 4, 393--407.
 
 \bibitem{F2} Floer, Andreas. The unregularized gradient flow of the symplectic action. Comm. Pure Appl. Math. 41 (1988), no. 6, 775--813.
 
 \bibitem{F3} Floer, Andreas. Morse theory for Lagrangian intersections. J. Differential Geom. 28 (1988), no. 3, 513--547. MR0965228.
 
 \bibitem{F4} Floer, Andreas. Cuplength estimates on Lagrangian intersections. Comm. Pure Appl. Math. 42 (1989), no. 4, 335--356.
 
 \bibitem{F5} Floer, Andreas. Witten's complex and infinite-dimensional Morse theory. J. Differential Geom. 30 (1989), no. 1, 207--221. MR1001276.
 
 \bibitem{F6} Floer, Andreas. Symplectic fixed points and holomorphic spheres. Comm. Math. Phys. 120 (1989), no. 4, 575--611.
 
 \bibitem{FM} Urs Frauenfelder, Agustin Moreno. \emph{On GIT quotients of the symplectic group, stability and bifurcations of symmetric orbits}. Preprint arXiv:2109.09147.
 
 \bibitem{FvK} Frauenfelder, Urs; van Koert, Otto. The restricted three-body problem and holomorphic curves. Pathways in Mathematics. Birkhäuser/Springer, Cham, 2018. xi+374 pp. ISBN: 978-3-319-72277-1; 978-3-319-72278-8.
 
 \bibitem{FvK2} Frauenfelder, Urs; van Koert, Otto. The Hörmander index of symmetric periodic orbits. Geom. Dedicata 168 (2014), 197--205.
 
\bibitem{Ginzburg} Ginzburg, Viktor L. The Conley conjecture. Ann. of Math. (2) 172 (2010), no. 2, 1127--1180.
 
\bibitem{henon} H\'enon, Michel. \emph{Numerical Exploration
of the Restricted Three-Body Problem. V. Hill's Case: Periodic
Orbits and Their Stability}, Astron.\,\&\,Astrophysics \textbf{1} (1969),
223--238.

\bibitem{HM} Howard, J. E. and MacKay, R. S. (1987). \emph{Linear stability of symplectic maps}, J. Math. Phys. 28, 1038-1051.
 
\bibitem{KABM} D.\ Koh, R.L.\ Anderson, I.\ Bermejo-Moreno. \emph{Cell-mapping orbit search for mission design at ocean worlds using parallel computing}. The Journal of the Astronautical Sciences 68 (1), (2021) 172-196.
 
\bibitem{KABM2018} D.\ Koh, R.L.\ Anderson, I.\ Bermejo-Moreno. \emph{Three-dimensional bifurcations in the circular restricted three-body problem}. AAS/AIAA Astrodynamics specialist conference, Snowbird, UT, August, 19–23 (2018) AAS 18–264. 

\bibitem[LTJ]{LTJ} Li Q, Tao Y, Jiang F. Orbital Stability and Invariant Manifolds on Distant Retrograde Orbits around Ganymede and Nearby Higher-Period Orbits. Aerospace. 2022; 9(8):454. https://doi.org/10.3390/aerospace9080454

\bibitem{kre1} Krein, M.: Generalization of certain investigations of A.M. Liapunov on linear
differential equations with periodic coefficients. Doklady Akad. Nauk USSR 73 (1950) 445-448.

\bibitem{kre2} Krein, M.: On the application of an algebraic proposition in the theory of monodromy matrices. Uspekhi Math. Nauk 6 (1951) 171-177.

\bibitem{K3} Krein, M.: On the theory of entire matrix-functions of exponential type. Ukrainian Math. Journal 3 (1951) 164-173.

\bibitem{K4} Krein, M.: On some maximum and minimum problems for characteristic numbers and Liapunov stability zones. Prikl. Math. Mekh. 15 (1951) 323-348.

\bibitem{Moser} Moser, J.: New aspects in the theory of stability of Hamiltonian systems. 
Comm. Pure Appl. Math. 11 (1958) 81-114. 
 
\bibitem{RS} J. Robbin, D. Salamon, \emph{The Maslov index for paths}, Topology 32, (1993), 827–844.

 
 \bibitem{W96} M.\,Wonenburger,
 \emph{Transformations which are products of two involutions}, J.\,Math.\,Mech.
 \textbf{16} (1996), 327--338.
 
\end{thebibliography}
\end{document}